\documentclass[12pt]{amsart}
\usepackage{amssymb}
\usepackage{amsfonts}
\usepackage{amsthm}
\usepackage{amscd}
\numberwithin{equation}{section}
\theoremstyle{plain}
\newtheorem{theorem}{Theorem}[section]
\newtheorem{proposition}[theorem]{Proposition}

\newtheorem{lemma}[theorem]{Lemma}
\theoremstyle{definition}
\newtheorem{definition}[theorem]{Definition}
\newtheorem{remark}[theorem]{Remark}
\theoremstyle{remark}
\newtheorem{claim}[theorem]{Claim}
\newcommand{\refE}[1]{(\ref{E:#1})}
\newcommand{\refS}[1]{Section~\ref{S:#1}}
\newcommand{\refSS}[1]{Section~\ref{SS:#1}}
\newcommand{\refT}[1]{Theorem~\ref{T:#1}}

\newcommand{\refP}[1]{Proposition~\ref{P:#1}}
\newcommand{\refD}[1]{Definition~\ref{D:#1}}

\newcommand{\refCl}[1]{Claim~\ref{C:#1}}
\newcommand{\refR}[1]{Remark~\ref{R:#1}}
\newcommand{\refL}[1]{Lemma~\ref{L:#1}}
\renewcommand{\a}{\ensuremath{\alpha}}
\newcommand{\g}{\ensuremath{\gamma}}
\newcommand{\w}{\ensuremath{\omega}}
\renewcommand{\b}{\ensuremath{\beta}}
\newcommand{\R}{\ensuremath{\mathbb{R}}}
\newcommand{\C}{\ensuremath{\mathbb{C}}}
\newcommand{\N}{\ensuremath{\mathbb{N}}}

\renewcommand{\P}{\ensuremath{\mathbb{P}}}
\newcommand{\Z}{\ensuremath{\mathbb{Z}}}
\newcommand{\K}{\ensuremath{\mathcal{K}}}

\newcommand{\cins}{\frac 1{2\pi\mathrm{i}}\int_{C_S}}

\newcommand{\cint}[1]{\frac 1{2\pi\mathrm{i}}\int_{#1}}
\newcommand{\cintl}[1]{\frac 1{24\pi\mathrm{i}}\int_{#1 }}

\newcommand{\A}{\mathcal{A}}
\newcommand{\hDT}{{h.d.t.}}
\renewcommand{\H}{\mathrm{H}}
\newcommand{\hL}{\mathrm{h.l.}}
\newcommand{\tr}{\mathrm{tr}}
\renewcommand{\i}{\mathrm{i}}
\renewcommand{\l}{\lambda}
\newcommand{\G}{{\mathcal{G}}}
\newcommand{\Gh}{\widehat{{\mathcal{G}}}}
\newcommand{\ga}{\mathfrak{g}}
\newcommand{\gb}{\overline{\mathfrak{g}}}
\newcommand{\gh}{\widehat{\mathfrak{g}}}

\newcommand{\KN} {Kri\-che\-ver-Novi\-kov}

\newcommand{\iK}{\ensuremath{1,\ldots, K}}
\newcommand{\Hl}[1][\lambda]{\mathcal{H}^{#1}}
\newcommand{\Fl}[1][\lambda]{\mathcal{F}^{#1}}
\newcommand{\Fn}[1][\lambda]{\mathcal{F}^{#1}}
\newcommand{\Fln}[1][n]{\mathcal{F}_{#1}^\lambda}

\newcommand{\Flns}[1][n]{F_{#1}^{\lambda,*}}

\newcommand{\Ah}{\widehat{\mathcal{A}}}
\newcommand{\La}{\mathcal{L}}
\renewcommand{\L}{\mathcal{L}}

\newcommand{\Lh}{\widehat{\mathcal{L}}}

\newcommand{\Dal}{\mathcal{D}_\lambda}
\newcommand{\Dalh}{\widehat{\mathcal{D}}_\lambda}
\newcommand{\Do}{\mathcal{D}^1}
\newcommand{\Dh}{\widehat{\mathcal{D}^1}}

\newcommand{\kndual}[2]{\langle #1,#2\rangle}
\newcommand{\piK}{P_1,\ldots, P_K}
\newcommand{\ord}{\operatorname{ord}}
\newcommand{\res}{\operatorname{res}}

\newcommand{\Cal}[1]{\mathcal{#1}}
\newcommand{\fpz}{\frac {d }{dz}}
\newcommand{\ldot}{\,.\,}
\newcommand{\dzl}{\,{dz}^\l}
\newcommand{\pfz}[1]{\frac {d#1}{dz}}
\newcommand{\de}{\delta}
\renewcommand{\d}{\delta}
\newcommand{\gli}{gl(\infty)}
\newcommand{\glih}{\widehat{gl}(\infty)}
\newcommand{\glib}{\overline{gl}(\infty)}
\renewcommand{\Re}{\mathrm{Re}}

\begin{document}
%\baselineskip=20pt

%\layout
%%%%%%%%%%%%%%%%%    private header  %%%%%%%%%%%%%%%%%%%%
\vspace*{-1cm}
\hbox{ }
{\hspace*{\fill} Mannheimer Manuskripte 265}

{\hspace*{\fill} math/0112116}

\vspace*{2cm}

\title[Local cocycles and central extensions]
{Local cocycles and central extensions for multi-point
algebras of Krichever-Novikov type}

\author[M. Schlichenmaier]
{Martin Schlichenmaier}
\address[Martin Schlichenmaier]{Department of Mathematics and 
  Computer Science, University of Mannheim, A5, 
         D-68131 Mannheim,
         Germany}
\email{schlichenmaier@math.uni-mannheim.de}
\begin{abstract}
Multi-point algebras of Krichever Novikov type for higher
genus Riemann surfaces are generalisations of the Virasoro
algebra and its related algebras.
Complete existence 
and uniqueness results for local 
2-cocycles defining almost-graded central extensions of the functions algebra,
the vector field algebra, and the differential operator algebra 
(of degree $\le 1$) are shown.
This is applied to the higher genus,  multi-point affine 
algebras to obtain uniqueness for almost-graded  central extensions of
the current algebra  of a simple
finite-dimensional Lie algebra. An earlier conjecture 
of the author 
concerning the central extension of the differential operator algebra
induced by  the 
semi-infinite wedge representations is proved.
\end{abstract}
\date{December 12, 2001}
\maketitle
%%%%%%%%%%%%%%%%%%%%%%%
% section 1
\section{Introduction}\label{S:intro}
%\input intro.tex
%%%%%%%
%Introduction
%%%%%%%%%%%%%%%%%%%%%%%%%%%%%%%%
%%%%%%%%%%%%%%%%%%%%%%%%%%%%%%%
Algebras of Krichever-Novikov type are important examples of 
infinite-dimensional associative algebras or Lie algebras. They generalize
the Witt algebra, its universal central extension (the Virasoro algebra) and
related algebras like the untwisted affine (Kac-Moody) algebra.
One way to describe the Witt algebra is to define it as the algebra of
those meromorphic vector fields on the Riemann sphere $S^2=\P^1(\C)$ 
which have only  poles at $0$ and $\infty$.
It admits a standard basis
$\ \{e_n=z^{n+1}\fpz,\ n\in\Z \}$.
The Lie structure is  the Lie bracket of vector fields.
One obtains immediately
\begin{equation*}
[e_n,e_m]=(m-n)e_{n+m}.
\end{equation*}
By introducing the degree $\deg(e_n):=n$ it becomes a graded 
Lie algebra.
In such a way all related algebras to the Witt algebra can be given as
 meromorphic  objects on $S^2$ which are holomorphic outside $0$ and $\infty$.
Of special importance besides the vector field algebra
are  the function algebra,
i.e. the algebra of Laurent polynomials $\C[z,z^{-1}]$,
the (Lie algebra) $\ga$-valued meromorphic functions,
i.e. the loop or current algebra $\ga\otimes\C[z,z^{-1}]$ with 
structure 
\begin{equation*}
[x\otimes z^n,y\otimes z^m]:=[x,y]\otimes z^{n+m}, \quad
x,y\in\ga,\quad n,m\in\Z,
\end{equation*}
and their central extensions, the Virasoro
algebra, the Heisenberg algebra and the untwisted
affine (Kac-Moody) algebras.

If one replaces $S^2$ by higher genus compact Riemann surfaces
(or equivalently by smooth projective curves over $\C$) and
allows, instead at two points, poles at a set $A$ of finitely many points
which is divided into two disjoint nonempty subsets $I$ and $O$,
one obtains in a similar way the algebras of Krichever-Novikov type
(see \refS{kn} for details).
For higher genus and two points the vector field algebra and the
function algebra were introduced by Krichever and Novikov
\cite{KNFa,KNFb,KNFc}, the corresponding affine algebras by Sheinman
\cite{Shea,Shaff}.
Its multi-point generalization
were given by the author 
\cite{Schlkn,Schleg,Schlce,SchlDiss},
including the function algebra, the vector field algebra, the
differential operator algebra, and the current algebra
\cite{Schlct,Schlhab,SchlShSu} and  their central 
extensions.

For the algebras related to the  Virasoro algebra the fact
that they are graded is of importance in many contexts. In particular,
this plays a role  in their representation theory (e.g. highest 
weight representations, Verma modules, etc.).
It turns out that a weaker concept, an almost-grading, will be
enough to guarantee the availability of certain methods in 
representation theory of infinite-dimensional algebras.
Almost-grading means that for pairs of homogeneous elements of
degree $n$ and $m$ the result is in a fixed range 
(not depending on $n$ and $m$) around the ``ideal'' value
$n+m$ (see \refD{almgrad}).
In the works cited above it is shown that there exists for any splitting
of $A$ into $I\cup O$ a grading, such that the algebras and
their modules are almost-graded. 

Krichever-Novikov type algebras
appear e.g. in string theory, in conformal field theory and
also in the theory of integrable models.
In particular, in closed string theory in the interpretation
 of the Riemann surface $M$  as possible world sheet of the theory,
the points in $I$ correspond to free incoming strings and the points
in $O$ to free outgoing strings.
The non-simply-connectedness of $M$ corresponds to string
creation, annihilation and interaction.
Furthermore these algebras have relations to moduli spaces,
e.g. \cite{SchlShWZ1,SchlShWZ2}.

In all the above-mentioned fields the passage to central extensions
of the algebras are of fundamental importance. Typically,
by some necessary regularization procedure one obtains only projective 
representations of the involved algebra which can be 
given as
linear representations of a suitable central extension.
Such a central extension is given by a 2-cocycle of
the Lie algebra cohomology with values in the trivial module.
For the representation theory it is fundamental to extend the almost-grading
to the central extension.
This requires that the defining 2-cocycle is local,
where we understand by a local cocycle a cocycle which 
vanishes if calculated for pairs of 
homogeneous elements of degree $n$ and $m$ if the sum $n+m$ lies outside
a certain fixed range (not depending on $n$ and $m$).

For the considered algebras there are certain cocycles
geometrically defined.
These cocycles are given in \refE{fung}, \refE{vecg}, and \refE{mixg}.
They are obtained by integration over cycles on the Riemann surface
with the points in $A$ removed.
If one chooses as integration cycle a cycle $C_S$ which separates the points
in $I$ from the points in $O$ one obtains a local cocycle
with respect to the almost-grading introduced by the splitting
$A=I\cup O$.

In this article I show that (up to coboundary) all local 
cocycles are scalar multiple of the above mentioned 
geometric cocycles obtained by integration along a separating
cycle.
The result is formulated in 
Theorems \ref{T:function}, \ref{T:vector},  
\ref{T:mix}, and \ref{T:diffop}.
In the function algebra case one obtains uniqueness only 
if one requires the cocycle to be a multiplicative or a
$\L$-invariant cocycle (see \refD{fmli}).
These properties are typically fulfilled in the applications under
consideration.
In particular, we obtain
\begin{equation*}
\dim \H^2_{loc,*}(\A,\C)=\dim \H^2_{loc}(\L,\C)=1,\quad
\dim \H^2_{loc}(\Do,\C)=3.
\end{equation*}
Here $\A$ denotes the function algebra, $\L$ the vector field
algebra, and $\Do$ the algebra of differential operators of degree
$\le 1$,
$\H_{loc}$ denotes the subspace of cohomology classes containing
at least one local cocycle, and $\H^2_{loc,*}(\A,\C)$ denotes local
cocycles which are (equivalently) either multiplicative or
$\L$-invariant.

Clearly, the 
classical case ($g=0$ and two points) 
is contained as a special case in the general results.
In the classical case the result for the Witt algebra is 
the well-known fact, that the Virasoro algebra is 
the universal central extension of the Witt algebra.
The statement for the differential operator algebra in the
classical case was shown by Arbarello, De~Concini, Kac and Procesi
\cite{ADKP}.
For the vector field algebra in higher genus with two points  
Krichever and Novikov supplied a proof of the uniqueness in 
a completely different manner than presented here.
Assuming that every cocycle is of geometric origin they
used ``discrete Baker-Akhieser functions'' to identify the
integration cycle \cite{KNFa,KNFb}.

The content of the article is as follows.
In \refS{kn} the  necessary basic informations about the geometric setup
and the studied algebras and its modules are given.
In \refS{geocyc}
central extensions and cocycles are studied.
In \refS{results} local cocycles are introduced and the main results
about uniqueness are formulated.
\refS{proofs} contains
the proofs.
The technique presented there involves the almost-grading and
consists essentially in setting-up a suitable recursion between
different levels.
In the vector field and differential operator algebra case the explicit
description of the basis elements via rational functions and 
theta function respectively is needed.

In \cite{SchlDiss,Schlct} the author formulated a conjecture
about the cocycle of the differential operator algebra associated to
a representation on the semi-infinite wedge forms of weight $\lambda$.
In \refS{wedge} it is shown that the conjecture follows from the
results obtained in this article (\refT{pullcyc}).
In particular, the identified cocycle extends to the whole
differential operator algebra of arbitrary degree.

\refS{affine} deals with an application to central extensions 
of current algebras $\ga\otimes \A$. In particular, if $\ga$ is a 
finite-dimensional simple Lie algebra 
any almost-graded central extension of  $\ga\otimes \A$
is obtained by 
a scalar multiple of a geometric cocycle 
for which the integration is over a 
separating cycle, see \refT{simple}.

There are some  articles addressing the different
question of 
determining the full cohomology space (or at least its
dimension)  of some of the
algebras considered here.
For the vector field algebra $\L$ see for example
results by Wagemann \cite{Wagkn,Wagdens} based on work
of Kawazumi \cite{KawGF}.
From these it follows
that $\dim \H^2(\L,\C)=2g+N-1$,
where $g$ is the genus of the Riemann surface $M$ and $N$ is the number of
points in $A$.
See also some earlier work of Millionshchikov
\cite{Mill} in which he proves finite-dimensionality.
Further there is the work of Getzler \cite{Getzdiff},
Wodzicki \cite{Woddiff},
and Li \cite{Li} on the differential operator algebra of all
degrees,
and Kassel and Loday \cite{Kaskd,KasLod}, Bremner \cite{Bremce,Brem4a},
and others for the current algebras.
These results can not be used in the theory of highest weight
representations of the algebras, because the
almost-grading (via the locality) cannot be  incorporated.
In general, the full cohomology spaces
are  higher dimensional.
Roughly speaking, the deRham cohomology of $M\setminus A$ 
is responsible for the 
Lie algebra cocycles.
For the classical case the deRham cohomology space is one-dimensional.
Hence, in this case (and only in this case) all 2-cocycle classes are
local classes and we recover the classical results.
But in general, to identify the local cocycle classes seems to be
a difficult task. The approach presented here is completely different.
We do not use the partial results on the general cohomology mentioned
above, but use a direct approach.
In addition, we deal systematically with even a  broader class
of algebras.

%\newpage
%%%%%%%%%%%%%%%%%%%%%%%%%%%%%%%%%%%
\section{The multi-point algebras of Krichever-Novikov type}\label{S:kn}
%\input kn.tex
%\section{The multi-point Krichever-Novikov algebras}\label{S:kn}
%\input kn.tex
%\section{The multi-point Krichever-Novikov algebras}\label{S:kn}

\subsection{Geometric set-up and the algebra structure}
$ $
\medskip

Let $M$ be a compact Riemann surface of genus $g$, or in terms
of algebraic geometry, a smooth projective curve over $\C$.
Let $N,K\in\N$ with $N\ge 2$ and $1\le K<N$. Fix
$\ 
I=(P_1,\ldots,P_K),$\ {and}\ $O=(Q_1,\ldots,Q_{N-K})$\ 
disjoint  ordered tuples of  distinct points (``marked points''
``punctures'') on the
curve. In particular, we assume $P_i\ne Q_j$ for every
pair $(i,j)$. The points in $I$ are
called the {\it in-points} the points in $O$ the {\it out-points}.
Sometimes we consider $I$ and $O$ simply as sets and set
$A=I\cup O$ as a set.

Let $\K$ be the canonical line bundle of $M$.
Its associated sheaf of local sections is the sheaf of
holomorphic differentials.
Following the common practice I will usually not
distinguish between a line bundle and its associated invertible sheaf
of section.
For every $\l\in\Z$ we consider the bundle
$\ \K^\l:=\K^{\otimes \l}$. Here we use the usual convention:
$\K^0=\Cal O$ is the trivial bundle, and $\K^{-1}=\K^*$
is the holomorphic tangent line bundle (resp.
the sheaf  of holomorphic vector fields).
After
fixing a theta characteristics, i.e. a bundle  $S$ with
$S^{\otimes 2}=\K$, we can also consider $\l\in \frac {1}{2}\Z$
with respect to the chosen theta characteristics.
In this article we will only need $\l\in\Z$.
Denote by $\Fl$ the (infinite-dimensional) vector space% 
\footnote{
For $\lambda=\frac 12+\Z$ we should denote the vector space by
$\mathcal{F}^\l_S$
and let $S$ go through all theta characteristics. 
}
of
global meromorphic sections  of $\K^\l$
 which are holomorphic
on $M\setminus A$.

Special cases, which are of particular interest to us, are
the quadratic differentials ($\l=2$),
the  differentials ($\l=1$),
the functions  ($\l=0$), and
the vector fields ($\l=-1$).
The space of functions I will also denote by $\A$ and the
space of vector fields by $\L$.
By  multiplying  sections with functions
we again obtain sections. In this way
the space $\A$ becomes an associative algebra and the spaces $\Fl$ become
$\A$-modules.

The vector fields in $\L$ operate on $\Fl$ by taking
the Lie derivative.
In local coordinates
\begin{equation}\label{E:Lder}
L_e(g)_|:=(e(z)\fpz)\ldot (g(z)\dzl):=
\left( e(z)\pfz g(z)+\l\, g(z)\pfz e(z)\right)\dzl \ .
\end{equation}
Here $e\in \L$ and $g\in \Fl$. To avoid cumbersome notation I
used the same symbol for the section and its representing
function.
If there is no danger of confusion I will do the same in the
following.

The space $\L$ becomes a Lie algebra with respect to \
the Lie derivative \refE{Lder}
and the spaces $\Fl$ become Lie modules over $\L$.
As usual I write $[e,f]$ for the bracket of the vector fields.
Its local form is
\begin{equation}\label{E:Lbrack}
[e(z)\fpz, f(z)\fpz]=
\left( e(z)\pfz f(z)- f(z)\pfz e(z)\right)\fpz \ .
\end{equation}

For the Riemann sphere ($g=0$) with quasi-global coordinate $z$
 and $I=(0)$ and $O=(\infty)$ the introduced 
function algebra is the algebra of Laurent polynomials 
$\C[z,z^{-1}]$ and the 
vector field algebra is
the Witt algebra, i.e.  the algebra whose universal central extension
is the Virasoro algebra.
We denote for short this situation as the
{\it classical situation}.
%%%%%%%%%%%%%%%%%%%%%%%%%%%%%%%%%%%%%%%%%%%%%%%%%

The vector field algebra $\L$ 
operates on the algebra $\A$ of functions
as derivations.
Hence it is possible to consider the semi-direct product
$\Do=\A\times \L$.
This Lie algebra is the algebra of differential operators of degree
$\le 1$ which are holomorphic on $M\setminus A$.
As vector space $\Do=\A\oplus\L$ and the Lie product is given as
\begin{equation}
[(g,e),(h,f)]:=(e\ldot h-f\ldot g,[e,f]).
\end{equation}
There is the short exact sequence of Lie algebras
\begin{equation}
\begin{CD}
0@>>>\A@>>>\Do@>>>\L@>>>0.
\end{CD}
\end{equation}
Obviously, $\L$ is also a  subalgebra of $\Do$.
The vector spaces $\Fl$ become $\Do$-modules by the canonical 
definition
\begin{equation}
(g,e)\ldot v=g\cdot v+e\ldot v,\quad v \in\Fl.
\end{equation}
By universal constructions algebras of differential operators of 
arbitrary degree can be considered 
\cite{SchlDiss, Schlct, Schlwed}.
There is another algebra of importance, the current algebra.
It will be defined in \refS{affine}.

%%%%%%%%%%%%%%%%%%%%%%%%%%%%%%%%%%%%%%

Let $\rho$ be a meromorphic differential which is holomorphic on 
$M\setminus A$ with exact pole order $1$ at the points in $A$ and given
positive residues at $I$ and given negative residues at $O$
(of course obeying the restriction 
$\sum_{P\in I}\res_P(\rho)+
\sum_{Q\in O}\res_Q(\rho)=0$) and purely imaginary periods.
There exists exactly one such $\rho$ (see \cite[p.116]{SchlRS}).
For $R\in M\setminus A$ a fixed point, the function
$u(P)=\Re\int_R^P\rho$ is a well-defined harmonic function.
The family of level lines
$ C_\tau:=\{p\in M\mid u(P)=\tau\},\ \tau\in\R $
defines a fibration of $M\setminus A$.
Each $C_\tau$ separates the points in $I$ from the points in $O$.
For $\tau\ll 0$ ($\tau\gg 0$) each level line $C_\tau$ is a disjoint union of
deformed circles $C_i$ around the points $P_i$, $i=\iK$ 
(of deformed circles $C_i^*$  around the points $Q_i$, $i=1,\ldots,N-K$).

For $f\in\Fl$ and $g\in\Fl[\mu]$ we have  $f\otimes g\in\Fl[\l+\mu]$.
In particular for $\mu=1-\l$ we obtain a meromorphic differential.
\begin{definition}\label{D:knpair}
The {\it Krichever-Novikov pairing} ({\it KN pairing}) is the
pairing between $\Fl$ and $\Fl[1-\l]$ given by
\begin{equation}\label{E:knpair}
\begin{gathered}
\Fl\times\Fn[1-\l]\ \to\ \C,
\\
\kndual {f}{g}:=\cint{{C_\tau}}f\otimes g
=\sum_{P\in I}\res_{P}(f\otimes g)=
-\sum_{Q\in O}\res_{Q}(f\otimes g),
\end{gathered}
\end{equation}
where $C_\tau$ is any non-singular level line.
\end{definition}
The last equality follows from the residue theorem.
Note that in \refE{knpair} the integral does not depend on
the level line chosen.
We will call any such level line or any cycle cohomologous to such
a level line a separating
cycle.
In particular, the KN pairing can be described as
\begin{equation}\label{E:knpairg}
\kndual {f}{g}=
\frac {1}{2\pi\i}\sum_{i=1}^K\int_{C_i}f\otimes g=
\cins f\otimes g=
-\frac {1}{2\pi\i}\sum_{i=1}^{N-K}\int_{C_i^*}f\otimes g.
\end{equation}

%%%%%%%%%%%%%%%%%%%%%%%%%%%%%%%%%%%%%%%%%%%%%%%%%%%%%%%%%%%%%
%%%%%%%%%%%%%%%%%%%%%%%%%%%%%%%%%%%%%%%%%%%%%%%%%%%%%%%%%%%%%%%%%%%%

\subsection{Almost-graded structure}$ $
\medskip

For infinite dimensional algebras 
and their representation theory a graded structure is usually
of importance to obtain structure results.
A typical example is given by the Witt algebra $W$. $W$ admits a
preferred set of basis elements given by
$\{e_n=z^{n+1}\fpz\mid n\in\Z\}$.
One calculates $[e_n,e_m]=(m-n)e_{n+m}$.
Hence $\deg(e_n):=n$ makes $W$ to a graded Lie algebra.

In our more general context the algebras will almost never be graded.
But it was observed by Krichever and Novikov 
in the two-point case that a weaker
concept, an almost-graded structure (they call it a quasi-graded 
structure), will be enough to develop an
interesting  theory of representations (highest weight representations, 
Verma modules, etc.).
\begin{definition}\label{D:almgrad}
(a) Let $\L$ be an (associative or Lie) algebra admitting a direct
decomposition as vector space $\ \L=\bigoplus_{n\in\Z} \L_n\ $.
The algebra $\L$ is called an {\it almost-graded}
algebra if (1) $\ \dim \L_n<\infty\ $ and (2)
there are constants $R$ and  $S$ with
\begin{equation}\label{E:eaga}
\L_n\cdot \L_m\quad \subseteq \bigoplus_{h=n+m+R}^{n+m+S} \L_h,
\qquad\forall n,m\in\Z\ .
\end{equation}
The elements of $\L_n$ are called {\it homogeneous  elements of degree $n$}.
\newline
(b) Let $\L$ be an almost-graded  (associative or Lie) algebra
and $\Cal M$ an $\L$-module with
$\ \Cal M=\bigoplus_{n\in\Z} \Cal M_n\ $
as vector space. The module $\Cal M$ is called an {\it almost-graded}
module, if
(1) $\ \dim \Cal M_n<\infty\ $, and
(2) there are constants $R'$ and  $S'$ with
\begin{equation}\label{E:egam}
 \L_m \cdot\Cal M_n\quad \subseteq \bigoplus_{h=n+m+R'}^{n+m+S'} \Cal M_h,
\qquad \forall n,m\in\Z\ .
\end{equation}
The elements of $\Cal M_n$ are called {\it homogeneous  elements of degree $n$}
\end{definition}
\noindent
By a {\it weak almost-grading} we understand an almost-grading without
requiring the finite-dimensionality of the homogeneous
subspaces.

For the 2-point situation for 
$M$ a higher genus Riemann surface and  $I=\{P\}$, $O=\{Q\}$
with $P,Q\in M$, Krichever and Novikov
introduced an almost-graded structure of the algebras and the modules
by exhibiting
special bases and defining their elements to be the
homogeneous elements.
In \cite{Schlce,SchlDiss} its multi-point
generalization was given, again  by
exhibiting a special basis.
(See also Sadov \cite{Sad} for some results in  similar directions.)

In more detail, 
for fixed $\l$ 
and for  every $n\in\Z$, and $i=1,\ldots,K$ a certain element
$f_{n,p}^\l\in\Fl$ is exhibited.
The $f_{n,p}^\l$ for $p=1,\ldots,K$ are a basis of a
subspace $\Fln$ and it is shown that
$$
\Fl=\bigoplus_{n\in\Z}\Fln\ .
$$
The subspace  $\Fln$ is called the {\it homogeneous subspace of degree $n$}.
%%%%%%%%%%%%%%%%%%%%%%%%%%%%%%%%%%%%%%%%%%%%%%%%%%%%%%%%%

The basis elements are chosen in such a way that they
fulfill the duality relation
with respect to the KN pairing \refE{knpair}
\begin{equation}\label{E:dual}
\kndual {f_{n,p}^\l} {f_{m,r}^{1-\l}}=
=\de_{-n}^{m}\cdot
\de_{p}^{r}\ . 
\end{equation}
This implies that the KN pairing is non-degenerate.

We will need as  additional 
information about the elements  $f_{n,p}^\l$ that
\begin{equation}\label{E:ordfn}
\ord_{P_i}(f_{n,p}^\l)=(n+1-\l)-\d_i^p,\quad i=1,\ldots,K .
\end{equation}
The recipe for choosing the order at the points 
in $O$ is such that up to a scalar multiplication
there is a unique such element which also fulfills \refE{dual}.
After choosing local coordinates $z_p$ at the points
$P_p$ the scalar can be fixed by requiring
\begin{equation}
{f_{n,p}^\l}_|(z_p)=z_p^{n-\l}(1+O(z_p))\left(dz_p\right)^\l\ .
\end{equation}

To give an impression of the type of conditions at $O$ let me consider
two cases. For
$N=K+1$ and $O=\{Q_{1}\}$ for  $g\ge 2$, $\l\ne 0,1$ and a generic
choice for the points in $A$  (or $g=0$ without any restriction)
we set
\begin{equation}\label{E:ordi}
\ord_{Q_1}(f_{n,p}^\l)=-K\cdot(n+1-\l)
+(2\l-1)(g-1)\ .
\end{equation}
For $N=2K$ and $O=(Q_1,Q_2,\ldots,Q_K)$
we set
\begin{equation}\label{E:ordin}
\begin{aligned}
\ord_{Q_i}(f_{n,p}^\l)=-(n+1-\l),\quad i=1,\ldots,K-1
\\
\ord_{Q_N}(f_{n,p}^\l)=-(n+1-\l)
+(2\l-1)(g-1).
\end{aligned}
\end{equation}
For $\l=0$ or $\l=1$ (and hence for all $\lambda$ in the case of genus $g=1$)
for small $n$ some modifications are  necessary. For $g\ge 2$
and for certain
values of $n$ and $\l$ such modifications are also needed if the 
points are not in generic positions.
See \cite{Schlce} for the general recipe.
By Riemann-Roch type arguments it is shown in \cite{Schlkn}
that there is up to a scalar multiple  only one such $f_{n,p}^\l$.

For the basis elements $f_{n,p}^\l$ in \cite{Schleg} explicit
descriptions in terms of rational functions (for $g=0$),
the Weierstra\ss\ $\sigma$-function (for $g=1$), and prime forms
and theta functions (for $g\ge1$) are given.
For a description using Weierstra\ss\ $\wp$-function, see
\cite{RDS}, \cite{SchlDeg}.
We will need  such a description at a certain
step in our proofs.

If $f\in\Fl$ is any element then it can be written as
$\ 
f=\sum_{m,r}'\a_{m,r}f_{m,r}^\l.
\ $
To simplify notation I will sometimes  use $\sum_{m,r}$ to denote
the double sum $\sum_{m\in\Z}\sum_{r=1}^K$.
The symbol $\sum'$ denotes that only finitely many terms
will appear in the sum.
Via \refE{dual} the coefficients can be 
calculated as
\begin{equation}
\a_{m,r}=\kndual {f}{f_{-m,r}^{1-\l}}=\cins f\otimes f_{-m,r}^{1-\l}.
\end{equation}
By considering the pole order at $I$ and $O$ a possible range for
non-vanishing $\a_{m,r}$ is given. 
A detailed analysis \cite{SchlDiss,Schlce} yields
\begin{theorem}\label{T:almgrad}
With respect to the above introduced grading the 
associative algebra
$\A$, and the Lie algebras 
$\L$, and $\Do$ are almost-graded and the modules $\Fl$ are
almost-graded modules over them.
In all cases the lower shifts in the degree of the
result (e.g. the numbers $R,R'$ in
\refE{eaga} and \refE{egam}) are zero.
\end{theorem}
The upper shifts can be explicitly calculated.
We will not need them here.
Let us abbreviate for terms of higher degrees as the one
under consideration the symbol $\hDT$.
By calculating the exact residues in the case of the lower bound
we obtain
\begin{proposition}\label{P:boundary}
%\begin{equation*}
\begin{alignat*}{2}
A_{n,p}\cdot A_{m,r}&=\de_p^r\cdot A_{n+m,r}+\hDT,\quad
&A_{n,p}\cdot f_{m,r}&=\de_p^r\cdot f_{n+m,r}+\hDT,\quad
\\
[e_{n,p}, e_{m,r}]&=\de_p^r\cdot(m-n)\cdot e_{n+m,r}+\hDT,\quad
&e_{n,p}\ldot f_{m,r}&=\de_p^r\cdot (m+\l n)\cdot f_{n+m,r}+\hDT.
\end{alignat*}
%\end{equation*}
\end{proposition}
%%%%%%%%%%%%%%%%%%%%%%%%%%%%%%%%%%%%%%%%%%%%%

Note that the grading does not depend on the numbering of the points in 
$I$.
Also the filtration $\mathcal{F}^\l_{(n)}$ introduced by 
the grading does not depend on renumbering the points in $O$ because
\begin{equation}
\mathcal{F}^\l_{(n)}:=\bigoplus_{m\ge n}\Fln[m]=
\{f\in\Fl\mid \ord_P(f)\ge n-\l,\ \forall P\in I\}.
\end{equation}
But this is an invariant description.
It also shows that a different recipe for the orders at $O$ will not
change the filtration.

\begin{remark}\label{R:inverted}
In the following we have also to consider the case when we interchange
the role played by $I$ and $O$. We obtain a 
different grading ${}^*$ introduced by  $I^*=O$.
This grading we call inverted grading.
For $N>2$ this not only a simple inversion and a translation.
Homogeneous elements of the original grading in general 
will not be homogeneous anymore and vice versa.
Denote the homogeneous objects and basis with  respect to the 
new grading also by ${}^*$. 
By considering the orders at the points $P_i$ and $Q_j$ 
and using \refE{dual} we obtain
\begin{equation}
\Fln \subseteq \bigoplus_{h=-\alpha n-L_1}^{-\alpha n+L_2}
\Flns ,\qquad
\Flns \subseteq \bigoplus_{h=-\beta n-L_3}^{-\beta n+L_4}
\Fln ,
\end{equation}
with $\alpha,\beta>0$ and $L_1,L_2,L_3,L_4$ numbers which do not
depend on $n$ and $m$.
\end{remark}

Let  me introduce the following
notation:
\begin{equation}\label{E:conc}
A_{n,p}:=f_{n,p}^0,\quad
e_{n,p}:=f_{n,p}^{-1},\quad
\w^{n,p}:=f_{-n,p}^1,\quad
\Omega^{n,p}:=f_{-n,p}^2 \ .
\end{equation}

%%%%%%%%%%%%%%%%%%%%%%%%%%%%%%%%%%
\section{Cocycles and central extensions}
\label{S:geocyc}
%\input geocyc.tex
%\section{Geometric cocycles}
%\label{S:geocyc}
%\section{Geometric cocycles}
In this section I consider central extensions of the above 
introduced algebras.
In quantum theory one is typically
forced (e.g. by regularization procedures) 
to consider projective representations of the algebras which
correspond to linear representations of centrally extended algebras.

Let $\G$ be any Lie algebra (over $\C$). A (one-dimensional)
central extension $\Gh$ is the middle term of a short exact sequence of
Lie algebras
\begin{equation}
\begin{CD}
0@>>>\C@>>>\Gh@>>>\G@>>>0
\end{CD}, 
\end{equation}
such that $\C$ is central in $\Gh$.
Two central extensions $\Gh_1$ and $\Gh_2$ are called equivalent if 
there is a Lie isomorphism $\varphi:\Gh_1\to\Gh_2$ such that
the diagram
\begin{equation}
\begin{CD}
0@>>>\C@>>>\Gh_1@>>>\G@>>>0
\\
@|@|@V{\varphi}VV@|@|
\\
0@>>>\C@>>>\Gh_2@>>>\G@>>>0
\end{CD}
\end{equation}
is commutative.

Central extensions are classified up to equivalence by
the second Lie algebra cohomology space 
$\H^2(\G,\C)$ (where $\C$ is considered as the
trivial module), i.e. by  2-cocycles up to coboundaries.
An antisymmetric map 
\begin{equation}
\g:\G\times \G\to \C
\end{equation}
is a 2-cocycle if
\begin{equation}\label{E:cocycle}
\g([f,g],h)+\g([g,h],f)+\g([h,f],g)=0,\quad\forall f,g,h\in\G.
\end{equation}
A  2-cocycle is a coboundary if there is a linear map 
$\phi:\G\to \C$ such that
\begin{equation}
\g(f,g)=\phi([f,g]),\quad \forall f,g\in\G.
\end{equation}
In the following we will only deal with 2-cocycles which we will just
call cocycles.
Given a cocycle $\g$ the central extension can be explicitly given 
by the vector space direct sum $\Gh:=\C\oplus \G$ with the Lie bracket given
by the  structure
equations (with $\widehat{e}:=(0,e)$ and $t:=(1,0)$)
\begin{equation}
[\widehat{e},\widehat{f}]:=\widehat{[e,f]}+\g(e,f)\cdot t,\quad
[t,\Gh]=0.
\end{equation}
In terms of short exact sequences we obtain
\begin{equation}
\begin{CD}
0@>>>\C@>{i_1}>>\Gh=\C\oplus\G@>{p_2}>>\G@>>>0.
\end{CD}
\end{equation}
Changing the cocycle by a coboundary corresponds to choosing a 
different linear lifting map of $p_2$ other than $i_2$.

In the following subsections we are considering 
cocycles for the algebras $\A$ (considered as
abelian Lie algebra), $\L$ and $\Do$ and the by the cocycles  defined central 
extensions.
In \refS{affine} we consider cocycles of the current algebras
(multi-point and higher genus).
For the classical situation the cocycles are either given 
purely algebraic in terms of structure constants or as
integrals (or residues) of objects expressed via the quasi-global 
coordinate $z$.
Typically they are not invariantly defined.
The classical expressions need some counter terms involving projective and 
affine connections.
\begin{definition}
Let $\ (U_\a,z_\a)_{\a\in J}\ $ be a covering of the Riemann surface
by holomorphic coordinates, with transition functions
$z_\b=f_{\b\a}(z_\a)$.
A system of local (holomorphic, meromorphic) functions 
$\ R=(R_\a(z_\a))\ $ resp. $\ T=(T_\a(z_\a))\ $
is called a (holomorphic, meromorphic) projective (resp. affine)
connection if it transforms as
\begin{equation}\label{E:pc}
R_\b(z_\b)\cdot (f_{\beta,\alpha}')^2=R_\a(z_\a)+S(f_{\beta,\alpha}),
\qquad\text{with}\quad
S(h)=\frac {h'''}{h'}-\frac 32\left(\frac {h''}{h'}\right)^2
\end{equation}
the Schwartzian derivative, respectively
\begin{equation}\label{E:ac}
T_\b(z_\b)\cdot f_{\beta,\alpha}'=T_\a(z_\a)+
\frac {f_{\beta,\alpha}''}{f_{\beta,\alpha}'}\ .
\end{equation}
Here ${}'$ denotes differentiation with respect to
the coordinate $z_\a$.
\end{definition}
It follows from \refE{pc} and \refE{ac} that  
the difference of two affine (projective) connections
is always a usual (quadratic) differential.
\begin{proposition}
Let $M$ be any compact Riemann surface.
\newline 
(a) There exists always
a holomorphic projective connection.
\newline
(b) Given a point $P$ on $M$ there exists always a meromorphic
 affine connection which is holomorphic
outside $P$ and has there at most a pole of order 1
\end{proposition}
\begin{proof}
(a) is a classical result, e.g. see \cite{Gun,HawSch}.
(b) is shown in \cite{SchlDiss,Schlwed}.
\end{proof}
For the following I will choose a fixed holomorphic projective 
connection $R^{(0)}$ and a fixed meromorphic affine connection 
$T^{(0)}$ with
at most a pole of order $1$ at the point $Q_1$.
All other connections with poles only at the points in $A$ can be
obtained by adding elements of $\Fl[1]$, resp. $\Fl[2]$, to these
reference connections.

%%%%%%%%%%%%%%%%%%%%%%%%%%%%%%%%%%%%%%%%%%%%%%
%%%%%%%%%%%%%%%%%%%%%%%%%%%%%%%%%%%%%%%%%%%
\subsection{Central extensions of the function algebra} $ $
\medskip

The function algebra considered as Lie algebra is abelian.
Hence any antisymmetric bilinear form will define a 2-cocycle.
For any $f,g\in \A$ and any linear form $\phi$ we obtain
$\phi([f,g])=0$. Hence, there will be no coboundary, i.e.
$\H^2(\A,\C)\cong \bigwedge^2\A$.

In the following we will consider cocycles which are of geometric 
origin.
Let $C$ be any differentiable cycle in $M\setminus A$ then 
\begin{equation}\label{E:fung}
\g:\A\times \A\to\C,\quad
\gamma_C(g,h):=\cint{C} gdh 
\end{equation}
is antisymmetric because
$0=\int_C d(gh)=\int_C gdh+\int_C hdg$.
Hence, this defines a cocycle.
Note that $C$ can be replaced by any
homologous cycle (assuming that it is  still a differentiable curve)
in $\H_1(M\setminus A,\Z)$,
because the differential $fdg$ is holomorphic on $M\setminus A$.
Any cocycle obtained via choosing a cycle $C$ in \refE{fung} is called
a {\it geometric cocycle}.
%%%%%%%%%%%%%%%%%%%%%%%%%%%%%%%%%%%%
\begin{definition}\label{D:fmli}
(a) A cocycle $\g$ for $\A$ is multiplicative
if it fulfills the ``cocycle condition'' for the associative algebra
$\A$, i.e.
\begin{equation}
\g(f\cdot g,h)+\g(g\cdot h,f)+\g(h\cdot f,g)=0,\quad \forall f,g,h\in\A.
\end{equation} 
(b)  A cocycle $\g$ for $\A$ is $\L$-invariant 
if 
\begin{equation}
\g(e.g,h)=\g(e.h,g),\quad \forall e\in\L,\ \forall g,h\in\A.
\end{equation}
\end{definition}
Both properties are of importance.
Below we will show that a cocycle of the function algebra
which is obtained via restriction from the differential
operator algebra will be  $\L$-invariant.
In \refS{wedge} we will show that cocycles obtained by
pulling back the standard cocycle of $\glib$
(see its definition there) 
via embeddings of $\A$ into $\glib$ respecting the almost-grading
will be multiplicative.
\newline
\begin{proposition}\label{P:geom}
The cocycle $\g_C$ \refE{fung} is multiplicative and  $\L$-invariant.
\end{proposition}
\begin{proof}
That $\g_C$ is multiplicative follows from $\int_C d(fgh)=0$ and Leibniz
rule.
To see the $\L$-invariance, first note that we have
$e\ldot dh=d(e\ldot h)$ for $e\in\L$ and $h\in\A$,
i.e. the differentiation and the Lie derivative commute.
Second, we have $e\ldot \omega=d(\omega(e))$ for  $e\in\L$ and
$\omega\in\Fl[1]$.
Both claims  can be directly
verified in local coordinates.
Now
\begin{equation*}
\int_C(e\ldot f)dg=\int_C e\ldot(fdg)-\int_C f\cdot(e\ldot dg)
=-\int_Cf\cdot(d(e\ldot g))=\int_C(e\ldot g)df.
\end{equation*}
In the first step we used 
$e.(a\otimes b)=(e\ldot a)\otimes b+a\otimes(e\ldot b)$
for $a\in\Fl$ and $b\in\Fl[\mu]$,
in the second step that  the first integral vanishes 
due to the fact that it is differential (using  $e\ldot \omega=d(\omega(e))$),
 and in the last step 
the antisymmetry of the cocycle.
\end{proof}
%%%%%%%%%%%%%%%%%%%%%%%%%%%%%%%%%%%%%%%%%%%%%%%%%%%%%%%%%%%%%%%
%%%%%%%%%%%%%%%%%%%%%%%%%%%%%%%%%%%%%%%%%%%%
\subsection{Central extensions of the vector field algebra} $ $
\medskip

In the classical situation there is up to equivalence and rescaling
only one nontrivial central extension of the Witt algebra, the
Virasoro algebra. In terms of generators $e_n$ the 
standard form of the cocycle is
\begin{equation}
\g(e_n,e_m)=\frac 1{12}(n^3-n)\de_{-m}^n.
\end{equation}
For the higher genus multi-point situation we consider for 
each cycle $C$ (or cycle class) with respect to the chosen
projective connection $R^{(0)}$
\begin{equation}\label{E:vecg}
\g_{C,R^{(0)}}(e,f):=\cintl{C} \left(\frac 12(e'''f-ef''')
-R^{(0)}\cdot(e'f-ef')\right)dz\ .
\end{equation}
Recall that we use the same letter for the vector field and its
local representing function.
This cocycle was introduced for the $N=2$ case by Krichever and
Novikov \cite{KNFa,KNFb}. 
As shown in \cite{SchlDiss,Schlce} it can be extended to the
multi-point situation.
There it was also shown that the integrand is indeed a differential and
that it defines a cocycle.

Next we consider coboundaries. A cocycle which is a coboundary can be given as
$D_\phi(e,f)=\phi([e,f])$ with a linear form $\phi$.
We have fixed a projective connection. If we choose another projective
connection $R$ which has only poles at $A$,
then $R=R^{(0)}+\Omega$ with a meromorphic quadratic differential with
poles only at $A$.
We calculate
\begin{equation}
\g_{C,R}(e,f)-
\g_{C,R^{(0)}}(e,f)=\cintl{C}\Omega(ef'-fe')dz=
\cintl{C}\Omega\otimes [e,f].
\end{equation}
This implies that the two cocycles are cohomologous.

The linear forms on $\L$ can be given in terms of the dual elements
of $e_{n,p}$. We can employ the KN pairing \refE{knpair} and can give 
$\phi$ by
\begin{equation}\label{E:exp}
\phi(e)=\kndual {W}{e}, \quad \text{with}\quad
W=\sum_{n\in\Z}\sum_r\beta_{n,r}\Omega^{n,r},
\end{equation}
see \refE{conc}.
Here the outer sum can reach indeed from $-\infty$ to $+\infty$.
Recall that 
\begin{equation}
\kndual {\Omega^{n,p}}{e_{m,r}}=\de_n^m\cdot \de_p^r,
\end{equation}
and that for a fixed $e$ only finitely many 
terms in \refE{exp} will be nonzero.
In this way we can give any coboundary by choosing such an 
infinite sum $W$.
Let us denote this coboundary by
\begin{equation}
D_W(e,f)=\kndual {W}{[e,f]}.
\end{equation}

We will call a cocycle a geometric cocycle if it can be represented 
as \refE{vecg} with a suitable cycle $C$ where the
reference connection might be replaced by an
meromorphic projective connection $R$.
\begin{remark}
One part of the Feigin-Novikov conjecture  says that every cocycle
of the vector field algebra is cohomologous to a linear combination
of geometric 
cocycles obtained by integration along the basis cycles in
$\H_1(M\setminus A,\Z)$ .
This (and the more general conjecture) was proven by 
Wagemann \cite{Wagkn,Wagdens} based on work of Kawazumi \cite{KawGF}.
We will not use this classification result in the following.
Instead we will show directly that every local cocycle (see \refD{local})
and more generally every cocycle  which is bounded from above
will be a geometric cocycle involving only the cycles
$C_1,C_2\ldots,C_K$.
\end{remark}
%%%%%%%%%%%%%%%%%%%%%%%%%%%%%%%%%%%%%%%%%%%%%%%%%%%%%%%%%%%%%%%%%%
%%%%%%%%%%%%%%%%%%%%%%%%%%%%%%%%%%%%%%%%%%%%%%%%%%%%%%%%%%%%%%%%%
\subsection{Central extensions of the differential operator algebra}$ $
\medskip

Due to the exact sequence of Lie algebras
\begin{equation}
\begin{CD}
0@>>>\A@>{i_1}>>\Do@>{p_2}>>\L@>>>0
\end{CD}
\end{equation}
every cocycle 
$\gamma^{(v)}$ of $\L$ will define via pull-back  a cocycle  
$p_2^*(\g)$ on $\Do$. 
Restricted to the subspace $\L$ in $\Do$ it will be exactly
the cocycle $\g^{(v)}$ and it will vanish if one of the arguments
is from $\A$.
We will denote this cocycle on $\Do$ also by $\g^{(v)}$.
 
The situation is slightly more complicated for the function algebra
$\A$ in $\Do$.
\begin{proposition}\label{P:linv}
A cocycle $\g^{(f)}$ of $\A$ can be extended to a cocycle of $\Do$ if and only
if $\g^{(f)}$ is $\L$-invariant
i.e.
\begin{equation}\label{E:linv}
\g^{(f)}(e.g,h)=\g^{(f)}(e.h,g),\quad \forall e\in \L,  \forall g,h\in\A. 
\end{equation}
\end{proposition}
\begin{proof}
Let $\tilde\g$ be a cocycle for 
$\Do$ and $\g$ its restriction to
$\A$. If we write down the cocycle condition for the elements
$e\in \L$ and $g,h\in\A$ we obtain \refE{linv}.
Vice versa: We define the extended bilinear map
\begin{equation}
\tilde\gamma:\Do\times\Do\to\C,\quad
\tilde\gamma((g,e),(h,f)):=\g^{(f)} (g,h).
\end{equation}
Clearly it is antisymmetric. We have to check the cocycle condition.
By linearity it is enough to do this for ``pure'' elements $(e,f,g)$.
If at least 2 of them are vector fields  or all of them are
functions then each of the terms in the cocycle relation vanishes
separately. It remains $e\in \L$ and $f,g\in\A$.
Because $[f,g]=0$ the cocycle condition is equivalent to \refE{linv}.
\end{proof}
By \refP{geom}
the geometric cocycles fulfill \refE{linv}.
Hence,
\begin{proposition}
The geometric cocycles 
$\g_C^{(f)}(f,g)=\cint{C}fdg$ can be
extended to $\Do$.
\end{proposition}
Let $\g$ be an arbitrary cocycle of $\Do$,
 and let $\g^{(f)}$ be its restriction
to $\A$ and $\g^{(v)}$ its restriction to $\L$ and both of them
extended by zero to $\Do$ again.
Then
$\g^{(m)}=\g-\g^{(f)}-\g^{(v)}$ will again be a cocycle.
It will only have nonzero values for $e\in\L$ and $f\in\A$ and fulfill
$\g^{(m)}(e,f)=-\g^{(m)}(f,e)$. We call $\g^{(m)}$ a mixing cocycle.
This decomposition of 
$\g=\g^{(f)}+\g^{(v)}+\g^{(m)}$ is unique.

Coboundaries for $\Do$ are given again by choosing  linear forms on $\Do$.
The dual spaces to the functions (vector fields) are given by 
the differentials (quadratic differentials) with the KN pairing
as duality.
Hence let
\begin{equation}\label{E:vw}
V=\sum_{n\in\Z}\sum_r\alpha_{n,r}\omega^{n,r},\quad
W=\sum_{n\in\Z}\sum_r\beta_{n,r}\Omega^{n,r}
\end{equation}
be possibly both-sided infinite sums then 
\begin{equation}
\phi((f,e))=\kndual {V}{f}+\kndual {W}{e}.
\end{equation}
The corresponding coboundary is given as
\begin{equation}
\phi([(g,e),(h,f)])=\phi((e\ldot h-f\ldot g,[e,f]))=
\kndual {V}{e\ldot h-f\ldot g}+
\kndual {W}{[e,f]}.
\end{equation}
This implies that the splitting into the three types remains if we pass to
cohomology.
The coboundary for $\g^{(v)}$ will be given by $W$, the
 coboundary for $\g^{(m)}$ will be given by $V$, and there is of course
no  coboundary for $\g^{(f)}$.

We want to study the mixing cocycles in more detail.
\begin{proposition}
Every bilinear form 
$
\g:\L\times \A\to\C
$
fulfilling 
\begin{equation}\label{E:mixco}
\g([e,f],g)-\g(e,f\ldot g)+\g(f,e\ldot g)=0\ , \forall e,f\in\L,\forall g\in\A
\end{equation}
defines by antisymmetric extension and by setting it zero on
$\A\times\A$ and $\L\times \L$  a mixing cocycle for $\Do$.
\end{proposition}
\begin{proof}
Let $\g$ be a bilinear form extended as described. Per construction
it is antisymmetric.
The only cocycle condition which does not trivially vanish is
the one involving two vector fields $e$ and $f$ and one function $g$. 
This cocycle condition is exactly \refE{mixco}.
\end{proof}

\begin{proposition}
Let $C$ be any cycle on the Riemann surface $M$.
And let $T^{(0)}$ be the meromorphic affine reference connection 
which has at most a pole of order 1 at $Q_{1}$  and 
is holomorphic elsewhere. Then
\begin{equation}\label{E:mixg}
\g_{C,T^{(0)}}(e,g)=-\g_{C,T^{(0)}}(g,e)=
\cint {C}\left(e\cdot g''+T^{(0)}\cdot (e\cdot g')\right)dz
\end{equation}
defines a mixing cocycle.
\end{proposition}
This has been shown in
\cite{SchlDiss} (see also \cite{Schlwed})
The addition of an affine connection is necessary because
otherwise the integrand would not be a differential.
As in the vector field case two cocycles obtained by 
different meromorphic affine connections with poles only at $A$
will be cohomologous.
Recall that $e\ldot g=e\cdot g'$, where the l.h.s. is
the Lie derivative with the vector field
and the r.h.s. is the multiplication with the local representing 
function.

As explained above the coboundaries can be given via 
$E_V(e,g)=\kndual {V}{e.g}$.
Again cocycles obtained via \refE{mixg}
with suitable affine connections are  called geometric cocycles.

In all three cases, of special importance are integration over the
cycles $C_1,C_2,\ldots, C_K$ around the points $P_i$, $i=1,\ldots,K$ and
integration over  the cycle $C_S=\sum_iC_i$.
The corresponding cocycles we will denote also by
$$
\g_i^{(f)},\  \g_{i,T}^{(m)},\ \g_{i,R}^{(v)},\quad i=1,\ldots, K,\
\quad\text{and}\quad
\g_S^{(f)},\ \g_{S,T}^{(m)},\ \g_{S,R}^{(v)}.
$$
The $S$ stands for the separating cycle $C_S$.
If the connection is the reference connection we will sometime
drop it in the notation.
A cocycle obtained via integration over a separating cycle I will call
a {\it separating cocycle}.

\begin{proposition}\label{P:indep}
In the following let $\g$ be  either the function cocycle \refE{fung},
the vector field cocycle \refE{vecg},
or the mixing cocycle \refE{mixg}.
\newline
(a) The cocycles $\g_i=\g_{C_i}$ for $i=1,\ldots,K$ are linearly independent
cohomology classes.
\newline
(b) The separating cocycle $\g_S$ is not cohomologous to zero.
\end{proposition}
\begin{proof}
The claim (b) follows from (a) because $\g_S=\sum_i\g_i$.
Now assume a linear relation 
$\sum_{i=1}^K\a_i[\g_i]=0$ in the cohomology space.
\newline
(i) We do first the function case. We evaluate this relation for the
pairs $(A_{-1,r},A_{1,r})$ with $r=1,\ldots,K$ and obtain $\a_r=0$
(there is no nontrivial coboundary). Hence, (a).
\newline
(ii) Mixing case:
The relation says there is a $V$ as in \refE{vw} such that 
$\sum_{i=1}^K\a_i\g_i=E_V$ (a possible $D_W$ will not contribute).
We evaluate this relation for pairs of elements $(e_{-n,r},A_{n,r})$
with $r=1,\ldots,K$ and obtain
\begin{multline}
\a_r\cdot n(n-1)-\kndual {V}{e_{-n,r}\ldot A_{n,r}}=
\a_r\cdot n(n-1)-
\kndual {V}{\sum_{h=0}^{L_2}\sum_t b_{(-n,r)(n,r)}^{(h,t)}A_{h,t}}
\\
=\a_r\cdot n(n-1)-B(e_{-n,r},A_{n,r})=0,\qquad
\end{multline}
with 
$B(e_{-n,r},A_{n,r}):=\sum_{h=0}^{L_2}\sum_t \a_{k,t}b_{(-n,r)(n,r)}^{(h,t)}$.
Here we used the almost-graded structure \refE{alf}
and the KN pairing \refE{knpair}.
\begin{claim}\label{C:indm}
$B(e_{-n,r},A_{n,r})=O(n)$.
\end{claim}
\noindent
Note that $L_2$ is a constant independent of $n$. Hence the
summation range will stay the same. But the coefficients may change with
$n$. We have to show that they are at most of order $n$ for $n\to\infty$.
This follows from the explicit description of the basis elements  of
$\Fl$
in terms of rational functions for $g=0$ and theta-functions and
prime forms for $g\ge 1$  given in \cite{Schleg}.
The details of the proof of the claim can be found in the  appendix.
Hence, $\a_r\cdot n(n-1)=O(n)$ which implies necessarily $\a_r=0$.
\newline
(iii) The vector field case is completely analogous with the modification that
as ``test pairs'' we take $(e_{-n,r},e_{n,r})$
and obtain 
\begin{equation}
\a_r\cdot (n+1)n(n-1)-C(e_{-n,r},e_{n,r})=0
\end{equation}
with \refE{almmix}. Again
\begin{claim}\label{C:indv}
$C(e_{-n,r},e_{n,r})=O(n)$.
\end{claim}
\noindent And we conclude as above.
\end{proof}

%%%%%%%%%%%%%%%%%%%%%%%%
\section{Uniqueness results for local cocycles}\label{S:results}
%\input results.tex
%%%%%%%%%%%%%%%%%%%%%%%%
%\section{Uniqueness results for local cocycles}\labe:{S:results}
%\input results.tex
%%%%%%%%%%%%%%%%%%%%%%%%%%%%%%%%%%
\begin{definition}\label{D:local}
(a) Let $\G=\bigoplus_{n\in\Z} \G_n$ be an almost-graded Lie algebra.
A cocycle $\g$ for $\G$ is called local it there exist $M_1,M_2\in\Z$ with
\begin{equation}\label{E:local}
\forall n,m\in\Z:\quad \g(\G_n,\G_m)\ne 0\implies 
M_2\le n+m\le M_1.
\end{equation}
(b) A cocycle  $\g$ for $\G$ is called bounded from above if there exists
 $M_1\in\Z$ with
\begin{equation}\label{E:bounded}
\forall n,m\in\Z:\quad \g(\G_n,\G_m)\ne 0\implies 
 n+m\le M_1.
\end{equation}
\end{definition}
If a cocycle is local the almost-grading of $\G$ can be extended to 
$\Gh=\C\oplus\G$
by defining $\deg \hat x=\deg x$ and $\deg t=0$.
Here $\hat x=(0,x)$ and $t=(1,0)$.
We call such an extension an almost-graded extension, or a local 
extension.
Krichever and Novikov \cite{KNFa} introduced the notion 
of local cocycles in the two point case and coined the name.
It might have been more suitable to use the name ``almost-graded cocycle''
instead of ``local cocycle''.
In any case, local cocycles are globally defined in contrast to their names.
%%%%%%%%%%%%%%%%%%%%%%%%%%%%%%%%%%%%%
\begin{theorem}\label{T:geomloc}
(a) The geometric cocycles 
$\g^{(f)}_S$, $\g^{(v)}_{S,R^{(0)}}$  and $\g^{(m)}_{S,T^{(0)}}$
are local cocycles which are bounded from above by zero.
\newline
(b) The geometric cocycles 
$\g^{(f)}_{C_i}$, $\g^{(v)}_{C_i,R^{(0)}}$  and $\g^{(m)}_{C_i,T^{(0)}}$
for $i=1,\ldots,K$ are bounded from above by zero.
\newline
(c) For an arbitrary meromorphic projective connection
$R$ and an arbitrary meromorphic affine connection $T$
which are holomorphic outside of $A$ the cocycles 
$\g^{(v)}_{S,R}$ and  $\g^{(m)}_{S,T}$
are local.
\end{theorem}
\begin{proof}
Recall that the index $S$ means integration over a separating cocycles.
The value of the above cocycles 
for homogeneous elements can be calculated by calculating residues 
at the points $\piK$.
Considering the order of the elements at these points we obtain 
that in case (b) the cocycles are
bounded from above by zero. 
Now $\g_S=\sum_i\g_i$, hence $\g_S$ is bounded from above by zero.
But equivalently the integration over a separating cycle can be
done by calculation of residues at the points $Q_1,\ldots,Q_{N-K}$.
This yields also a lower bound for them.
(See \cite{SchlDiss} for explicit formulas for the lower bounds).
Hence (a) follows.
As long as we add meromorphic 1-differentials (resp. 
quadratic differentials)  which have only poles at the points  of $A$ 
to the affine (resp. projective)
reference connection 
the bounds for the cocycles will change but they
will stay local.
The upper bound zero will not  change if we add only
1-differentials (resp. quadratic differentials) with maximal pole order 1
(resp. pole order 2) at the points in $I$.
\end{proof}

We call a cohomology class a local cohomology class if it contains 
a cocycle which is local.
This implies that  by choosing a suitable
lift of the elements of $\G$ to  $\Gh$ the almost-grading
of $\G$ can be extended to $\Gh$.
If $\g_1$ and $\g_2$ are local then the sum $\g_1+\g_2$ will also be
local.
Hence the local cohomology classes will be a subspace of $\H^2(\G,\C)$ which
we denote by   $\H^2_{loc}(\G,\C)$.
Note that not necessarily every element in a local cohomology class will
be a local cocycle.
%%%%%%%%%%%%%%%%%%%%%%%%%

\newpage

\begin{theorem}\label{T:function}
(a) A cocycle $\g$ for the function algebra $\A$ which is
either multiplicative or
$\L$-invariant is local if and only if it is  a multiple of the 
separating cocycle, i.e. there exists $\a\in\C$ such that
\begin{equation}
\g(f,g)=\alpha\g_S(f,g)=\frac {\alpha}{2\pi\mathrm{i}}\int_{C_S} fdg.
\end{equation}
(b) A local cocycle  will be bounded by zero and for  the values at the upper
bound we have
\begin{equation}\label{E:funzero}
\g(A_{-n,r},A_{n,s})=\alpha\cdot n\cdot\delta_r^s,
\quad\text{with}\quad
\alpha=\g(A_{-1,r},A_{1,r})
\end{equation}
for any $r=1,\ldots,K$.
\end{theorem}
As a consequence we obtain immediately
\begin{theorem}\label{T:fml}
(a) A local cocycle for the function algebra which is multiplicative
is also $\L$-invariant and vice versa.
\newline
(b) Denote by $\H^2_{loc,*}(\A,\C)$ the subspace of
cocycles which are local and multiplicative (or equivalently
local and differential), then 
 $\dim \H^2_{loc,*}(\A,\C)=1$.
\end{theorem}
Note that for the algebra $\A$ there are no nontrivial coboundaries.

%%%%%%%%%%%%%%%%%%%%%%%%%%%%%%%%%%%%%%%%%%%%%%%%%%%%%%%%%

\begin{theorem}\label{T:vector}
(a) A  cocycle for the vector field algebra 
$\g$ is a local cocycle if and only if
$\gamma$ is the sum of a multiple of the
separating cocycle with projective connection $R^{(0)}$ and of
a coboundary $D_W$, i.e. there exist $\a\in\C$
and  $W=\sum_{n=M_1}^{M_2}\sum_r\b_{n,r}\Omega^{n,r}$ such that
\begin{equation}
%\begin{gathered}
\g(e,f)=\alpha\g_{S,R^{(0)}}(e,f)+D_W(e,f),
\quad
%\\
\text{with}\quad D_W(e,f)=\kndual {W}{[e,f]}.
%\end{gathered}
\end{equation}
(b) If $\alpha\ne 0$ then 
$\ 
\g(e,f)=\alpha\g_{S,R}(e,f), 
\ $
with a projective connection  $R$ which has only poles at the points
in $A$.
\newline
(c) If $\g$ is a local cocycle which is
bounded from above by zero  then at  level zero
the cocycle is given by
\begin{equation}\label{E:veczero}
\begin{gathered}
\g(e_{n,r},e_{-n,s})=\left(\frac{(n+1)n(n-1)}{12}\cdot\alpha+
n b_{r}\right)\delta_s^r\ .
\\
\text{with}\quad
\a:=2\g(e_{2,r},e_{-2,r})-4\g(e_{1,r},e_{-1,r}),\ \text{and}\quad 
b_{r}:=\g(e_{1,r},e_{-1,r}).
\end{gathered}
\end{equation}
Here $\a$ can be calculated with respect to any $r$.
\newline
(d) $\dim \H_{loc}^2(\L,\C)=1$.
\end{theorem}

%%%%%%%%%%%%%%%%%%%%%%%%%%%%%%%%%%%%%%%%%%%%%%%%%%%%%%%%
\begin{theorem}\label{T:mix}
(a) A  mixing cocycle for the differential operator algebra 
$\g$ is a local cocycle if and only if
$\gamma$ is the sum of a multiple of the
separating cocycle with affine connection $T^{(0)}$ and of
a coboundary $E_V$, i.e. there exist $\a\in\C$ and  
$V=\sum_{n=M_1}^{M_2}\sum_r\b_{n,r}\omega^{n,r}$ such that
\begin{equation}
\g(e,g)=\alpha\g_{S,T^{(0)}}(e,g)+E_V(e,g),
\quad
\text{with}\quad E_V(e,g)=\kndual {V}{e.g}.
\end{equation}
(b) If $\alpha\ne 0$ then 
$\ \g(e,g)=\alpha\g_{S,T}(e,g) \ $
with an affine connection  $T$ which has only poles at the points
in $A$.
\newline
(c) If $\g$ is a local cocycle which is
bounded from above by zero then at  level zero the
cocycle is given by
\begin{equation}\label{E:mixzero}
\begin{gathered}
\g(e_{-n,r},A_{n,s})=\left(n(n-1)\alpha+n\cdot b_{r}\right)\cdot\delta_s^r\ .
\\
\a:=1/2\left(\g(e_{1,r},A_{-1,r})+\g(e_{-1,r},A_{1,r})\right),
\quad\text{and}\quad
b_{r}:=\g(e_{-1,r},A_{1,r}).
\end{gathered}
\end{equation}
Here $\a$ can be calculated with respect to any $r$.
\newline
(d) The subspace of local cohomology classes which 
are given by  mixing cocycles is one-dimensional.
\end{theorem}
%%%%%%%%%%%%%%%%%%%%%%%%%%%%%%%%%%%%%

\refT{function}, \refT{vector}, and \refT{mix} will be proved in the
following section.
In addition some more statements about cocycles which are
bounded from above are given. It will turn out
(Theorems \ref{T:fbound}, \ref{T:mbound}, \ref{T:vbound})  that
they are geometric cocycles involving as integration paths 
only the cycles $C_i$  around the points $\piK$.
This implies that for $K=1$ this path will be a separating cycle.
Hence,
\begin{proposition}
Let $K=1$ (e.g. $N=2$). 
\newline
(a) Every multiplicative or $\L$-invariant cocycle for the 
function algebra which is bounded from above will be local.
\newline
(b) Every cocycle for $\L$ or $\Do$ which is bounded from above is
cohomologous to a local cocycle.
\end{proposition}
\noindent
By passing to the inverted grading (see \refR{inverted}) the 
proposition is also true if we replace ``bounded from above'' by
bounded from below.

An arbitrary cocycle for the differential operator algebra can
uniquely be  decomposed
into 3 cocycles of fixed type.
Hence, we obtain as a corollary of the above theorems
\begin{theorem}\label{T:diffop}
A cocycle $\gamma$ for the differential operator algebra 
is local if and only if it is a linear combination of
the  cocycle obtained by
extension of the separating 
cocycle for the function algebra, the cocycle obtained 
by pulling back the  separating vector field cocycle and the 
separating mixing  cocycle with  meromorphic
affine and projective connections $T$ and $R$
holomorphic outside  $A$
\begin{equation}
\g\ =\ a_1\g^{(f)}_S+a_2\g^{(m)}_{S,T}+ 
a_3\g^{(v)}_{S,R} +E_V+D_W,
\end{equation}
and coboundary terms $E_V+D_W$.
If $a_2,a_3\ne 0$ then the coboundary terms can be
included into the connections $R$ and $T$.

(b) The subspace of $\H^2_{loc}(\Do,\C)$ of cocycles cohomologous to
local cocycles 
is three-dimensional
and is generated by the cohomology classes
of the separating cocycles of function, mixing and vector field type.
\end{theorem}

In the classical case, i.e. $M=\P^1(\C)$ and $A=\{0,\infty\}$, 
the statement about the vector field algebra is the 
well-known fact that the Witt algebra possess a one-dimensional universal
central extension, the Virasoro algebra.
In the standard description the cocycle is local.
In this case for the differential operator algebra the result was 
proved by Arbarello, De Concini, Kac and Procesi \cite{ADKP}.

For the higher genus  two-point situation  
the result for the vector field algebra was proved by 
Krichever and Novikov \cite{KNFa,KNFb}
starting from  the general assumption that the
cocycle will be a geometric cocycle and determining the
defining cycle using ``discrete Baker-Akhieser functions''.
The method of the proof is completely different.
I will not use in my proof their results.
Indeed, I will obtain an independent proof of it.

%%%%%%%%%%%%%%%%%%%%%%%%%%%%%%%%%%
\newpage
\section{The proofs}
\label{S:proofs}
\subsection{Multiplicative cocycles for the function algebra}
\label{SS:function}
$ $
\medskip

In the following  $\g$ denotes a multiplicative
cocycle for the function algebra which is bounded from above.
For a pair $(A_{n,p},A_{m,r})$ we call the sum $l=n+m$ the {\it level}
of the pair.
The pairs of level $l$ can be written as
$\g(A_{n,p},A_{-n+l,r})$.

We   make descending recursion on the level $l$.
First we will show that starting at a level $l>0$ for which the 
values of the cocycle will be zero for  $l$ and all higher
levels the values will  also be zero for all levels between 1 and  $l$.
Then we will show that for all levels less than zero the cocycle values are
determined by its values at level 0.
Finally, we analyse the level zero.
In particular it will turn out that all possible values
for level zero can be realized  by suitable linear combinations
of the  geometric cocycles $\g_r$, $r=1,\ldots,K$.
We conclude that $\g$ itself is a linear combination.
Boundedness from below will only allow  a combination for which all
coefficients are the same.
   
\begin{lemma}\label{L:rs}
The elements $\g(A_{n,r},A_{-n+l,s})$ of level $l$ for $r\ne s$ are 
universal linear combinations 
of  elements of level $\ge (l+1)$.
\end{lemma}
\begin{proof}
By the multiplicativity
$$
\g(A_{0,r}\cdot A_{n,r},A_{-n+l,s})+
\g(A_{n,r}\cdot A_{-n+l,s},A_{0,r})+
\g(A_{-n+l,s}\cdot A_{0,r},A_{n,r})=0.
$$
We replace the products with the help of the almost-grading for $\A$, 
i.e by
\begin{equation}\label{E:fgradc}
A_{n,r}\cdot A_{m,s}=\delta_r^s \cdot A_{n+m,r}+
\sum_{h=n+m+1}^{n+m+L}\sum_t a_{(n,r),(m,s)}^{(h,t)} A_{h,t},
\end{equation}
where $a_{(n,r),(m,s)}^{(h,t)}\in\C$, and $L$ is the upper bound
for the almost grading. As usual any 
summation range 
over the second index is
$\{1,\ldots,K\}$.
Hence  for $r\ne s$
$$
\g(A_{n,r}+\hDT,A_{-n+l,s})+
\g(\hDT,A_{0,r})+
\g(\hDT,A_{n,r})=0\ .
$$
Here $\hDT$ should denote linear combinations of 
elements of degree which do not contribute to the levels under
considerations.
This implies $\g(A_{n,r},A_{-n+l,s})$ 
can be expressed as linear combinations of 
values of the cocycle  of higher level than  $l$.
The coefficients appearing in this linear combination only
depend on the the geometric situation, i.e. on the 
structure constants of the algebra  and not on the
the cocycle under consideration.
This should be understood by the term ``universal linear
combination'' in the theorem.
\end{proof}

\begin{remark}\label{R:universal}
In the following we will  use the phrase
``can be expressed by elements of higher level'',
``determined by higher level'', or simply 
``$=\hL$'' to denote that it is a universal linear combination
of cocycle values for pairs of homogeneous 
elements of level higher than the level under consideration.
By the level we understand the sum of the degree of the two
arguments.
In particular, if two cocycles 
are given by higher level and they coincide in  higher levels, they will
coincide also for the elements under consideration.
\end{remark}
\begin{lemma}\label{L:0l}
The value $\g(A_{0,r},A_{l,r})$ can be
 expressed by elements of level $\ge l+1$.
\end{lemma}
\begin{proof}
By the multiplicativity
$$
\g(A_{0,r}\cdot A_{0,r},A_{l,r})+
\g(A_{0,r}\cdot A_{l,r},A_{0,r})+
\g(A_{l,r}\cdot A_{0,r},A_{0,r})=0\ .
$$
Using the almost-grading we obtain
$$
\g(A_{0,r},A_{l,r})+2\cdot\g(A_{l,r},A_{0,r})=\hL\ .
$$
By the antisymmetry of the cocycle the claim follows.
\end{proof}
We  do not need it in the following. But for completeness let me
note
\begin{lemma}\label{L:funlemma}
$
\g(1,f)=0,\quad\forall f\in\A\ .
$
\end{lemma}
\begin{proof}
From
$
\g(1\cdot 1,f)+
\g(1\cdot f,1)+\g(f\cdot 1,1)=0
$ we conclude 
$0=\g(1,f)+2\g(f,1)=\g(f,1)$.
\end{proof}

By \refL{rs} only the case $r=s$ is of importance at the level $l$.
Hence, to simplify notation we will suppress in the following
the second index.
Starting from
$$
\g(A_k\cdot A_n,A_m)+
\g(A_n\cdot A_m,A_k)+
\g(A_m\cdot A_k,A_n)=0
$$
we obtain
\begin{equation}\label{E:motherf}
\g(A_{k+n},A_{m})+\g(A_{n+m},A_{k})+\g(A_{m+k},A_n)=\hL\ .
\end{equation}
We specialize this for $m=-1$ and $m=1$:
\begin{gather}
\g(A_{k+n},A_{-1})+\g(A_{n-1},A_k)+\g(A_{k-1},A_n)=\hL,\label{E:mm1}
\\
\g(A_{k+n},A_1)+\g(A_{n+1},A_k)+\g(A_{k+1},A_n)=\hL,\label{E:mp1}
\end{gather}
and set in \refE{mm1} $k=l-n+1$ and in 
\refE{mp1} $k=l-n-1$ ($l$ denotes the level) to obtain
\begin{gather}
\g(A_{l+1},A_{-1})+\g(A_{n-1},A_{l-(n-1)})+\g(A_{l-n},A_n)=\hL,\label{E:mm1n}
\\
\g(A_{l-1},A_1)+\g(A_{n+1},A_{l-(n+1)})+\g(A_{l-n},A_n)=\hL\ .\label{E:mp1n}
\end{gather}
Subtracting \refE{mm1n} from \refE{mp1n} we obtain the 
recursion formula
\begin{equation}\label{E:recfun}
\g(A_{n+1},A_{l-(n+1)})=\g(A_{n-1},A_{l-(n-1)})-\g(A_{-1},A_{l+1})
+\g(A_1,A_{l-1})+\hL\ .
\end{equation}
If we set  $n=-m$ and $k=l$ in \refE{motherf} we obtain
$$
\g(A_{l-m},A_m)+\g(A_0,A_l)+\g(A_{l+m},A_{-m})=\hL.
$$
From \refL{0l} it follows that $\g(A_0,A_l)$ is  of higher level, hence
\begin{equation}
\g(A_m,A_{l-m})=-\g(A_{-m},A_{l+m})+\hL\ .
\end{equation}
For $m=1$ we obtain $\g(A_{1},A_{l-1})=-\g(A_{-1},A_{l+1})+\hL$,
 which we can plug into
\refE{recfun} to obtain
\begin{equation}
\g(A_{n+1},A_{l-(n+1)})=\g(A_{n-1},A_{l-(n-1)})+2\g(A_1,A_{ l-1})+\hL.
\end{equation}
Hence, the knowledge of $\g(A_0,A_l)$ and $\g(A_1,A_{l-1})$
 will fix the complete
cocycle at level $l$ by the knowledge of the higher levels.
But $\g(A_0,A_l)$ itself is fixed by higher level (\refL{0l}), hence 
$\g(A_1,A_{l-1})$, or equivalently $\g(A_{-1},A_{l+1})$ will fix everything.

First we consider the level  $l=0$ and obtain the recursion
$$
\g(A_{n+1},A_{-(n+1)})=\g(A_{n-1},A_{-(n-1)})+2\g(A_1,A_{-1})+\hL.
$$
This implies
\begin{equation}\label{E:funz1}
\g(A_n,A_{-n})=n\cdot\g(A_1,A_{-1})+\hL.
\end{equation} 
%\newpage
\begin{lemma}\label{L:fundet}
The level $l$ for $l\ne0$ is completely determined by higher levels.
\end{lemma}
\begin{proof}
First consider $l>0$.
We have to show that $\g(A_1,A_{l-1})$ is determined by higher levels.
For $l=1$ we obtain $\g(A_1,A_{l-1})=\g(A_1,A_0)$
 which is determined by higher 
level (see \refL{0l}). For $l=2$ we obtain $\g(A_1,A_{l-1})=\g(A_1,A_1)=0$ 
by the 
antisymmetry. Hence we can assume $l>2$.
We set in \refE{motherf} $k=l-r-1,n=1,m=r$ and obtain 
\begin{equation}\label{E:runfun}
\g(A_1,A_{l-1})+\g(A_r,A_{l-r})-\g(A_{r+1},A_{l-r-1})=\hL\ .
\end{equation} 
Set $m:= \frac {l-2}{2}$ for $l$ even or
$m:= \frac {l-1}{2}$ for $l$ odd.
We let $r$ run trough $1,2,\ldots,m$ and obtain from 
\refE{runfun} $m$ equations.
The first equation will always be
$$
2\cdot \g(A_1,A_{l-1})-\g(A_2,A_{l-2})=\hL.
$$
The last equation will depend on the parity of $l$.
For $l$ even and $r=m$ the last term 
on the l.h.s. of \refE{runfun}
will be $\g(A_{\frac l2},A_{\frac l2})$, which vanishes.
For $l$ odd the last term of the last equation will
coincide with the second term. Hence
$$
\g(A_1,A_{l-1})+2\cdot \g(A_{\frac {l-1}{2}},A_{\frac {l+1}{2}})=\hL
$$
will be the last equation.
In this case we divide it  by 2.
All these equations are added up. As result we obtain
$$
(m+\epsilon)\cdot \g(A_1,A_{l-1})=\hL,\
$$
where $\epsilon$ is 1 for $l$ even $1/2$ for $l$ odd.
This shows the claim for $l>0$.

For $l<0$ note that we can equally determine $\g(A_{-1},A_{l+1})$ to fix the 
cocycle.
Now the arguments work completely in the same way as above.
The claim for $l=-1,-2$ follows immediately.
We plug $k=l-r+1,n=-1,m=r$ into \refE{motherf} and obtain
\begin{equation}\label{E:runfunm}
\g(A_{-1},A_{l+1})+\g(A_r,A_{l-r})-\g(A_{r-1},A_{l-r+1})=\hL.
\end{equation}
We set $m:=\frac {-l-2}{2}$ for $l$ even and 
$m:=\frac {-l-1}{2}$ for $l$ odd and consider the equation
 \refE{runfunm} for $r=-1,-2,\ldots, -m$.
They have the similar structure as for $l>0$ and we can add them
up again to obtain the statement about $\g(A_{-1},A_{l+1})$.
\end{proof}
\begin{proposition}\label{P:bounddet}
Let $\gamma$ be multiplicative cocycle which is bounded
 from above then: 
\newline
(a) $0$ is also 
an upper bound, i.e 
$\g(\A_n,\A_m)=0$ for $n+m>0$.
\newline
(b) It is  determined by its value at level $0$.
\newline
(c) 
The level zero  is given as  
\begin{equation}\label{E:funz}
\g(A_{n,r},A_{-n,s})=n\cdot \delta_s^r\cdot \alpha_r,\quad
\text{with}\ \a_r:= \g(A_{1,r},A_{-1,r})
\end{equation}
 for $\quad n\in\Z,\ r,s=1,\ldots, K$.
\end{proposition}
\begin{proof}
Assume $\g$ to be bounded. If $M>0$ is an upper bound then by
\refL{fundet} and \refL{rs} its values at the level $M$ are
linear combinations of levels $>M$. Hence they are also vanishing
on level $M$ and finally $0$ is also an upper bound. This proves (a).
Part (b) follows again from the above lemmas. 
For $r\ne s$ the Equation \refE{funz} follows from \refL{rs}.
For $r=s$ this is \refE{funz1} which has to be applied for each $r$ 
separately.
\end{proof}
%%%%%%%%%%%%%%%%%%%%%%%%%%%%%%%%%%%%%%%%%%%%%%%%%%%%%%%%%%%%%%%%%%%%%
\begin{theorem}\label{T:fbound}
The space of multiplicative cocycles for the function algebra which are
bounded from above is $k-$dimensional.
A basis is given by the cocycles 
\begin{equation}
\g_{i}(f,g)=\cint{C_i}f dg,\quad i=1,\ldots, K.
\end{equation}
\end{theorem}
\begin{proof}
From the  \refP{bounddet}
it follows that the space is at most $K$-dimensional.
By \refT{geomloc} and \refP{geom} the 
geometric cocycles
$\gamma_i$ are bounded from above and 
multiplicative. Hence they are elements of this space.
By \refP{indep} they are linearly independent, hence a basis.
\end{proof}
%%%%%%%%%%%%%%%%%%%%%%%%%%%%%%%%%%%%
\begin{proof}[Proof of \refT{function} (multiplicative case)]
If $\gamma$ is a multiple of the separating cocycle then it is local
(see \refT{geomloc}).
Now assume that $\gamma$ is  multiplicative and local. Hence it is bounded from
above and can be written as $\g=\sum_{i}\alpha_i\g_i$. 
We have to show that $\a_1=\a_2=\cdots=\a_K$.
For $K=1$ claim (a) is immediate.
If we interchange the role of $I$ and $O$ we obtain 
the inverted grading 
(see \refR{inverted}) which we denote by ${}^*$. 
A cocycle which is bounded from 
below with respect to the old grading is bounded from above
with respect to the new grading.
Denote by $C_i^*$ circles around the points $Q_i$ for $i=1,\ldots, N-K$,
and by $\g_i^*(f,g)=\cint {C_i^*}fdg$ the corresponding geometric cocycle.
Using \refT{fbound}
we obtain $\g=\sum_{i=1}^{N-K}\alpha_i^*\g_i^*$ with certain
$\alpha_i^*\in\C$.
Again if $N-K=1$ then there is just one cocycle which is then
a separating cocycle and (a) is proven.
Hence we assume $N-K>1$.
By subtracting the two presentation of the same cocycle and
regrouping the summation we obtain
\begin{equation}\label{E:funsum}
0=(\a_1+\a_1^*)\sum_{i=1}^K\g_i +\sum_{k=2}^K (\a_k-\a_1)\g_k 
-\sum_{k=2}^{N-K} (\a_k^*-\a_1^*)\g_k^*.
\end{equation} 
Here we used $\sum_{i=1}^K\g_i+\sum_{i=1}^{N-K}\g_i^*=0$.

For each $k=2,\ldots, K$ separately we take the pair 
of functions $f_n$ and $g_n$ which are uniquely defined  for 
infinitely many $n$ with $n\gg 0$ by 
the conditions 
$\ord_{P_k}(f_n)=-n$, $\ord_{P_1}(f_n)=n-g$
and 
$\ord_{P_k}(g_n)=n$, $\ord_{P_1}(f_n)=-n-g$,
the requirement that they are holomorphic elsewhere and that  with
respect to the chosen local coordinate $z_k$ at $P_k$ the leading
coefficient is 1.
Then $\sum_{k=1}^K \g_k(f_n,g_n)=0$ because the
elements to not have poles at the points $Q_j$.
All terms in the sum \refE{funsum} are zero with the exception of
$\g_k(f_n,g_n)=n$.
This implies that $\a_k=\a_1$.
In a completely analogous  way $\a_k^*=\a_1^*$.
Hence it remains a multiple of the  separating 
cocycle. This can only be zero if $\a_1=-\a_1^*$.
This shows (a).
The explicit form \refE{funzero} 
of the level zero follows from \refP{bounddet} and
the fact that for the geometric cocycle we have
$\g_S(A_{-1,r},A_{1,r})=1$.
Hence (b).
\end{proof}
%%%%%%%%%%%%%%%%%%%%%%%%%%%%%%%%%%%%%%%%%%%%%%%%%%%%%%%%%%%

\subsection{$\L$-invariant  cocycles for the function algebra}
$ $
\medskip

In this subsection I consider cocycles of the
function algebra which are $\L$-invariant, i.e.
\begin{equation}\label{E:linvn}
\g(e.g,h)-\g(e.h,g)=0,\quad  e\in\L, g,h\in\A.  
\end{equation}
Cocycles which are obtained via restriction of cocycles of the
algebra $\Do$ of differential operators are of this type,
see \refP{linv}.
Also the geometric cocycles have this property (\refP{geom}).

By the almost-graded action of $\L$ on $\A$ we have
\begin{equation}\label{E:alf}
e_{n,r}.A_{m,s}=\delta_r^s \cdot m\cdot A_{n+m,r}
+\sum_{h=n+m+1}^{n+m+L_2}\sum_t 
 b_{(n,r),(m,s)}^{(h,t)} A_{h,t}
\end{equation}
with $b_{(n,r),(m,s)}^{(h,t)}\in\C$, and $L_2$ the upper bound
for the almost-graded structure.

Using \refE{linvn} we get
$$
\g_(e_{n,r}.A_{m,s},A_{p,t})=
\g_(e_{n,r}.A_{p,t},A_{m,s}),
$$
and with the 
almost-graded structure
\begin{equation}\label{E:diffmother}
\delta_s^r\cdot  m \cdot\g(A_{m+n,r},A_{p,t})=\delta_r^t\cdot   p\cdot 
 \g(A_{p+n,r},A_{m,s})+\hL.
\end{equation}
For $r=t\ne s$ we obtain
$p\cdot \g(A_{p+n,r},A_{m,s})=\hL$ for any $p$. This implies that
\refL{rs} is also true for 
$\L$-invariant 
cocycles.
Hence it is enough to consider $s=r=t$. We will drop again the
second index and obtain
\begin{equation}\label{E:diffrec}
m\cdot \g(A_{m+n},A_p)=p\cdot \g(A_{p+n},A_m)+\hL.
\end{equation}
If we set $n=0$ we obtain
\begin{equation}
(m+p)\cdot\g(A_m,A_p)=\hL.
\end{equation}
Note that $m+p$ is the level. Hence for level $l\ne 0$ everything is
determined by higher levels. This is 
\refL{fundet} now for $\L$-invariant cocycles.

Let us assume that $\g$ is bounded from above then (as above) it will 
also be  bounded by zero.
For level 0
we set $p=-(n+1)$ and $m=1$  in \refE{diffrec} and obtain
with the antisymmetry of the cocycle 
\begin{equation}
\g(A_{n+1},A_{-(n+1)})=(n+1)\cdot\g(A_{1},A_{-1}),
\end{equation}
which corresponds to \refE{funz1}.
The proofs of \refP{bounddet}, of \refT{fbound}
and of \refT{function} rely only
on these lemmas and the relation \refE{funz1}.
Hence we obtain that they are also valid if we replace ``multiplicative''
by ``$\L$-invariant''.
In particular, we obtain 
\begin{proposition}\label{P:boundl}
(a) The space of $\L$-invariant  cocycles for the function algebra which are
bounded from above is $K$-dimensional.
A basis is given by the cocycles 
\begin{equation}
\g_{i}(f,g)=\cint{C_i}f dg,\quad i=1,\ldots, K.
\end{equation}
(b) A bounded cocycle for the function algebra is multiplicative
if and only if it is $\L$-invariant.
\end{proposition}

%%%%%%%%%%%%%%%%%%%%%%%%%%%%%%%%%%%%%%%%%%%%%%%%%%%%%%%%%%%%%%
%\input mixing.tex
\subsection{Mixing local cocycles for the 
differential operator algebra}
\label{SS:mixing}
$ $
\medskip

In this subsection I consider those cocycles defined for the
differential operator algebra 
$\Do$ which  vanish on the subalgebra
$\A$ of functions and the subalgebra  $\La$ of vector fields.
We start from the cocycle relation for $e,g\in\L$ and $g\in\A$
\begin{equation}
\g([e,f],g)-\g(e,f.g)+\g(f,e.g)=0,
\end{equation}
which we evaluate for the basis elements
\begin{equation}\label{E:basmix}
\g([e_{k,r},e_{n,s}],A_{m,t})-\g(e_{k,r},e_{n,s}. A_{m,t})
+\g(e_{n,s},e_{k,r}.A_{m,t})=0\ .
\end{equation}
We use the almost-graded structure \refE{alf} and 
\begin{equation}\label{E:almmix}
[e_{k,r},e_{n,s}]=\delta_r^s \cdot (n-k)\cdot e_{k+n,r}+
+\sum_{h=n+m+1}^{n+m+L_3}\sum_t 
 c_{(n,r),(m,s)}^{(h,t)} e_{h,t}
\end{equation}
with $c_{(n,r),(m,s)}^{(h,t)}\in\C$, and $L_3$ the upper bound
for the almost-graded structure of $\L$.
Again we want to make induction on the level of the elements.
If we plug  \refE{alf} and \refE{almmix} into \refE{basmix} we obtain as
relation on level $l=(n+m+k)$ 
\begin{equation}\label{E:basemix}
\delta_r^s \cdot (n-k)\cdot\g(e_{k+n,r},A_{m,t})-
\delta_s^t\cdot m\cdot\g(e_{k,r},A_{m+n,s})+
\delta_r^t \cdot m\cdot\g(e_{n,s},A_{m+k,t})=\hL.
\end{equation}
If all $r,s,t$ are mutually distinct, this does not produce any
relation on this level.
For $s=t\ne r$, $m=-1$ and $n=p+1$
we obtain  $\g(e_{k,r},A_{p,s})=\hL$.
Hence,
\begin{lemma}\label{L:mixmix}
The cocycle values 
$\g(e_{k,r},A_{p,s})$ for $r\ne s$ can be expressed as 
universal linear 
combinations of cocycle values of higher level.
\end{lemma}
Again we use the phrase ``universal linear combination''
to denote the situation explained in \refR{universal}.

This shows that it 
is enough to consider elements with the same second index.
We will drop it in the notation.
The equation \refE{basemix} can now be written as:
\begin{equation}\label{E:motherm}
(n-k)\cdot \g(e_{k+n},A_{m})-
m\cdot \g(e_{k},A_{m+n})+
m\cdot \g(e_{n},A_{m+k})=\hL.
\end{equation}
Setting $k=0$ in \refE{motherm}
yields
\begin{equation}\label{E:motbas}
(n+m)\g(e_n,A_m)=m\g(e_0,A_{m+n})+\hL\ .
\end{equation}
\begin{lemma}
(a) If the level $l=(n+m)\ne 0$ then
\begin{equation}
\g(e_n,A_m)=\frac {m}{n+m}\cdot \g(e_0,A_{n+m})+\hL\ .
\end{equation}
(b) $\g(e_n,A_0)$ for all $n\in\Z$  is given 
by higher levels.
\end{lemma}
\begin{proof}
Part (a) is obtained by dividing \refE{motbas} 
by $n+m\ne 0$.
For $n\ne 0$ we obtain  (b) by setting $m=0$.
We set $m=1$ and $n=-1$ in   \refE{motbas}  and get (b) also for
$n=0$.
\end{proof} 
Hence, as long as the level $l\ne 0$, there is for each level only one
parameter which can be adjusted, then everything is
fixed by the higher levels.

It remains to deal with the level 0 case.
We set  $k=-n-m$ in \refE{motherm} and obtain
\begin{equation}
(2n+m)\cdot \g(e_{-m},A_m)-m\cdot \g(e_{-(n+m)},A_{n+m})+
m\cdot \g(e_n,A_{-n})=\hL.
\end{equation} 
We specialize further to $m=1$ and $m=-1$ 
\begin{gather}\label{E:mix1}
(2n+1)\g(e_{-1},A_{1})-\g(e_{-(n+1)},A_{n+1})+
\g(e_n,A_{-n})=\hL
\\ \label{E:mix2}
(2n-1)\g(e_{1},A_{-1})+\g(e_{-(n-1)},A_{n-1})
-\g(e_n,A_{-n})=\hL.
\end{gather}
Adding these equations yields
\begin{equation}
\g(e_{-(n+1)},A_{n+1})=
\g(e_{-(n-1)},A_{n-1})+(2n-1)\cdot \g(e_1,A_{-1})+
(2n+1)\cdot \g(e_{-1},A_{1})+\hL.
\end{equation}
Recall that $\g(e_0,A_0)=\hL$ hence the values on level zero
are uniquely fixed by $\g(e_1,A_{-1})$ and  $\g(e_{-1},A_{1})$.
We use 
\begin{equation}
\alpha:=1/2\left(\g(e_1,A_{-1})+\g(e_{-1},A_{1})\right),
\qquad
\beta_0:=-\g(e_{-1},A_{1}),
\end{equation}
and obtain
\begin{equation}
\g(e_{-(n+1)},A_{n+1})=\g(e_{-(n-1)},A_{n-1})+
2(2n-1)\alpha-2\beta_0\ .
\end{equation}
The starting elements of the recursion are
\begin{equation}
\g(e_{0},A_{0})=\hL,\qquad
\g(e_{-1},A_{1})=-\beta_0.
\end{equation}
By induction this implies
\begin{equation}
\g(e_{-n},A_{n})=n(n-1)\alpha-n\beta_0+\hL.
\end{equation}
\begin{proposition}\label{P:zerom}
(a) For a mixing cocycle the level 0 is fixed by the data
\begin{equation}
\alpha_r:=1/2 \left(\g(e_{1,r},A_{-1,r})+\g(e_{-1,r},A_{1,r})\right),\quad
\beta_{0,r}=-\g(e_{-1,r},A_{1,r}),
\end{equation}
for $r=\iK$ via
\begin{equation}\label{E:zerom}
\g(e_{-n,r},A_{n,s})=(n(n-1)\alpha_r-n\beta_{0,r})\cdot\delta_s^r+\hL,
\end{equation}
where $\hL$ denotes a universal linear combination of cocycle values
of level $> 0$.

(b) A mixing cocycle which is bounded from above is uniquely given 
by the collection of values
\begin{equation}\label{E:valmix}
\g(e_{1,r},A_{-1,r}),\g(e_{-1,r},A_{1,r}), \ r=1,\ldots, K,
\quad
\g(e_{0,r},A_{n,r}),\ n\in\Z\setminus\{0\}.
\end{equation}
\end{proposition}
\begin{proof}
Our analysis above works for every $r$ separately. 
\refL{mixmix} gives the statement for $r\ne s$.
This shows (a).

Let $\g_1$ and $\g_2$ be two cocycles bounded from above,
which have the same set of values \refE{valmix}.
Let $L$ be a common upper bound.
Recall that at level $l\ne 0$ the elements of this level are fixed 
as certain universal linear combinations  of elements of level
higher than $l$ and the element $\g(e_0,A_l)$.
Hence, the two cocycles coincide for every level $l> 0$.
For level 0 we use  part (a), hence they coincide on level 0
and further on on every level.
\end{proof}
%%%%%%%%%%%%%%%%%%%%%%%%%%%%%%%%%%%%%%%%%
\begin{theorem}\label{T:mbound}
(a) Let $\g$ be a mixing cocycle which is bounded from above by $M$
then there exist $\a_r\in \C$, $r=1,\ldots,K$  and  
a formal sum of 1-differentials
\begin{equation}
V=\sum_{n=-\infty}^{n=M}\sum_s\beta_{n,s}\w^{n,s},\quad \b_{n,s}\in \C,
\end{equation}
such that
\begin{equation}
\g=\sum_{r=1}^K\a_r\g_r+E_V,\qquad
E_V(e,g)=\kndual{V}{e.g}\ .
\end{equation}
(b) The cohomology space of mixing cocycles bounded from above
is $K$-dimensional and generated by the classes
$[\g_r],\ r=1,\ldots,K$.
\end{theorem}
Recall that $\g_r$ denotes the cocycle obtained by   
\refE{mixg} where we integrate over $C_r$ using our
fixed reference affine connection $T^{(0)}$.
\begin{proof}
In view of \refP{zerom} for proving (a) it is enough to show that we can realize by 
such a combination all values 
\begin{equation}
\g(e_{0,r},A_{n,r}), n\in\Z, n\le M, n\ne 0,\quad
\g(e_{1,r},A_{-1,r}),\ 
\g(e_{-1,r},A_{1,r}),\ r=1,\ldots, K.
\end{equation} 
Consider  a cocycle given as such a linear combination. 
We will show that we can recursively adjust all
parameter to realize all possible values.
The affine connection $T^{(0)}$ is fixed.
It does not have poles at the points in $I$.
Recall the orders of the basis elements
\begin{equation}
\ord_{P_r}(e_{n,r})=n+1,\quad
\ord_{P_r}(A_{n,r})=n,\quad
\ord_{P_r}(w^{n,r})=-n-1.
\end{equation}
The orders increase by 1 at the points $P_s$ with $s\ne r$.
The highest level is $M$. 
Assume $M>0$. We set all $\b_{m,r}=0$ for $m>M$.
The first set of values we have to realize are $\g(e_{0,r},A_{M,r})$
But in this case the first term in the expression does not
contribute, only
$E_V(e_{0,r},A_{M,r})$ will contribute.
But this term calculates to
\begin{equation}
\sum_{n=-\infty}^M\sum_s\b_{n,s}\kndual {\w^{n,s}} {e_{0,r}\ldot A_{M,r}}
= \sum_{n=-\infty}^M\sum_s\b_{n,s}\kndual{\w^{n,s}}
{M\cdot A_{M,r}+\hDT}
=\beta_{M,r}\cdot M.
\end{equation}
As long as $M>0$ we can divide by $M$ and  obtain the 
prescribed value $\g(e_{0,r},A_{M,r})$ for any $r$.
For $M-1>0$ we can do the same and obtain
\begin{equation}
\g(e_{0,r},A_{M-1,r})=(M-1)\beta_{M-1,r}+
\sum_{s=1}^K \beta_{M,s}\cdot b_{(0,r),(M-1,r)}^{(M,s)},
\end{equation}
where the coefficients $b_{(0,r),(M-1,r)}^{(M,s)}$
are the structure constants introduced in \refE{alf}.
This can be done recursively  as long as 
the level is $>0$.

For the level 0 we 
have also a contribution of the first term.
We calculate
\begin{equation}
\g(e_{1,r},A_{-1,r})=2\alpha_r+\beta_{0,r}+\hDT,
\quad
\g(e_{-1,r},A_{1,r})=-\beta_{0,r}+\hDT,
\end{equation}
where the higher degree terms are already determined.
This implies, that by setting
\begin{equation}
\a_r=1/2(\g(e_{1,r},A_{-1,r})+\g(e_{-1,r},A_{1,r}))+\hDT,
\quad
\beta_{0,r}=-\g(e_{-1,r},A_{1,r})+\hDT
\end{equation}
the level zero will have the prescribed values.

For $l<0$ the argument to determine
$\beta_{l,r}$ will work as for $l>0$ with the 
only modification, that we pick-up
additional  elements due to the first integral (which also involves the
expansion of $T^{(0)})$.
\newline
(b) follows now from (a) by \refP{indep}.
\end{proof}
%%%%%%%%%%%%%%%%%%%%%%%%%
\begin{proof}[Proof of \refT{mix}.]
First recall that a separating cocycle with finite sum in the coboundary
is local (see \refT{geomloc}).
To prove the opposite 
we use the same technique as in the proof of
\refT{function} presented in \refSS{function}.
The arguments are completely similar up to the formulation of
the equivalent of Equation \refE{funsum}.
The corresponding equation is
\begin{equation}\label{E:mixsum}
0=(\a_1+\a_1^*)\sum_{i=1}^K\g_i +\sum_{k=2}^K (\a_k-\a_1)\g_k 
-\sum_{k=2}^{N-K} (\a_k^*-\a_1^*)\g_k^*+E_V-E_{V^*}.
\end{equation}
Here 
\begin{equation}
V=\sum_{-\infty}^{n=M_1}\sum_r\beta_{n,r}\w^{n,r},\quad
V^*=\sum_{-\infty}^{n=M_2}\sum_r\beta_{n,r}^*\w_*^{n,r},
\end{equation}
and ${}^*$ denotes the opposite grading.
Again if $K=1$ (or $N-K=1$) 
the integration cycle will be the separating cycle.

Similar as there we consider for infinitely many $n$ with $n\gg 0$ 
a function $g_n$ and a vector field $e_n$
defined by 
$\ord_{P_k}(g_n)=n$, $\ord_{P_1}(g_n)=-n-g$,
$\ord_{P_k}(e_n)=-n+1$, $\ord_{P_1}(e_n)=n-3g$, and 
$\ord_{Q_1}(e_n)=1$ 
(this is due to the possible pole of $T^{(0)}$ at $Q_1$).
which are holomorphic elsewhere and which have with
respect to the chosen local coordinate $z_k$ at $P_k$ the leading
coefficients 1.
Now $\g_S(e_n,g_n)=0$ because there are no poles at the points in $O$.
Furthermore $\g_k(e_n,g_n)=n(n-1)$.
Hence,
\begin{gather}\label{E:asymp}
(\a_k-\a_1)n(n-1)+A(e_n,g_n)=0,\quad\text{with}
\\
A(e_n,g_n):=\cins\left(
\sum_{m=-(4g+1)}^{M_1}\sum_t\beta_{m,t}\w^{m,t}(e_ng_n')-
\sum_{m=0}^{M_2}\sum_t\beta_{m,t}^*\w^{m,t}_*(e_ng_n')\right).
\end{gather}
\begin{claim}\label{C:claimmix}
$A(e_n,g_n)=O(n)$.
\end{claim}
First note that independent of $n$ only the same finite 
summation ranges appear.
The lower bounds follows from 
$\ord_{P_1}(e_ng_n')=(-4g-1)$.
We have to show that the coefficients of the Taylor series expansion
of $e_ng_n'$ at the points of $A$ is of at most order $n$ if $n\to\infty$.
This follows from the explicit description of the basis elements  of
$\Fl$
in terms of rational functions for $g=0$ and theta-functions and
prime forms for $g\ge 1$  given in \cite{Schleg}.
The details of the proof of the claim can be found in the  appendix.
 
We  divide 
\refE{asymp}
by $n$ and obtain $(\a_k-\a_1)(n-1)+O(1)=0$ for $n\gg 0$, which 
implies $\a_k=a_1$.
This can be done for $k=2,\ldots, K$ and in the same manner for the
pair of points $(Q_1,Q_k)$.
Hence, $\gamma=\alpha_1\sum_k\g_k+E_V$.
Now $\g$ and $\g_S=\sum_k\g_r$  are local hence $E_V$ is a 
coboundary which is local. This implies that the sum $V$ is finite
(see \refP{finiteness} below).
This shows (a). If $\alpha=\alpha_1\ne 0 $ then we can write
$T:=T^{(0)}+(1/\a) V$ and obtain $\g=\alpha\cdot\g_{S,T}$. This shows (b).
Part (c) is only a specialization of \refP{zerom} and
the fact that for the separating cocycle 
$\g_S(e_{-1,r},A_{1,r})=0$  and $\g_S(e_{1,r},A_{-1,r})=2$. Part (d) follows
from (a) using \refP{indep}.
\end{proof}
\begin{proposition}\label{P:finiteness}
(a) Let  $\g_S^{(m)}$ be the separating mixing cocycle.
Then $\g=\g_S^{(m)}+E_V$ with 
$V=\sum_{-\infty}^{m=M}\sum_r\a_{m,r}\omega^{m,r}$
is local if and only if $V$ is a finite sum.
\newline
(a) Let  $\g_S^{(v)}$ be the separating vector field cocycle.
Then $\g=\g_S^{(v)}+D_W$ with 
$W=\sum_{-\infty}^{m=M}\sum_r\b_{m,r}\Omega^{m,r}$
is local if and only if $W$ is a finite sum.
\end{proposition}
\begin{proof}
I will only proof (a). The proof of (b) is completely analogous.
That finiteness of the sum implies locality follows from \refT{geomloc}.
For the opposite direction assume locality of $\g$.
With $\g$ also $\g-\g_S^{(m)}=E_V$ will be local.
Hence it is enough to proof the claim for $\g=E_V$
with 
$V=\sum_{-\infty}^{m=M}\sum_r \a_{m,r}\omega^{m,r}$.
We might even assume that $M$ is suitable negative.
Assume that $V$ is not finite.
Let $m_0$ be such that (1) $E_V(e_{n,r},A_{m,s})=0$ for
$(n+m)\le m_0$, (2) there exists an $r$ with $\a_{m_0,r}\ne 0$ and
(3) $m_0<0$.
We use the condition   $E_V(e_{k,r},A_{m_0-k,r})=0$ for 
all $k\ge 0$ and calculate with the almost graded structure \refE{alf}
\begin{equation}\label{E:indeg}
\a_{m_0,r}+\frac {1}{m_0-k}\sum_{h=m_0+1}^{m_0+L_2}
\sum_t\a_{h,t}a_{(k,r)(m_0-k,r)}^{(h,t)}=0,\quad k\in\N_0.
\end{equation}
If $L_2=0$ (which is the case for the classical situation) this
already implies that $\a_{m_0,r}=0$ in contradiction to the
assumption. Now assume $L_2>0$ then \refE{indeg} 
gives a homogeneous system of infinitely many independent equations for the 
$\a_{h,t}$ ($m_0+1\le h\le m_0+L$, $1\le t\le K$)  and  $a_{m_0,r}$.
We obtain only the trivial solution $\a_{m_0,r}$ in contradiction 
to the assumption.
\end{proof}

%%%%%%%%%%%%%%%%%%%%%%%%%%%%%%%%%%%%%%%%%%%%%%%%%%%%%%%%%%%%%%
%\input vector.tex
\subsection{Cocycles for the 
vector field algebra}
\label{SS:vector}
$ $
\medskip

For the cocycles of the vector field algebra the 
statements and the proofs are quite similar
to the mixing cocycle case.
Instead of  affine connections  projective connections will
appear.

If we plug the almost-grading \refE{almmix} into 
the cocycle condition 
\begin{equation}\label{E:vecbas}
\g([e_{n,r},e_{m,s}],e_{p,t})+
\g([e_{m,s},e_{p,t}],e_{n,r})+
\g([e_{p,t},e_{n,r}],e_{m,s})=0
\end{equation}
for triples of basis elements
we obtain
\begin{equation}
\delta^s_r\cdot (m-n)\g(e_{n+m,r},e_{p,t})+
\delta^s_t\cdot (p-m)\g(e_{m+p,s},e_{n,r})+
\delta^t_r\cdot (n-p)\g(e_{n+p,r},e_{m,s})=\hL.
\end{equation}
Again, if $s,t,r$ are mutually distinct this does not produce any
relation on level $n+m+p$.

For $s=r\ne t$ we obtain
\begin{equation}
(m-n)\g(e_{n+m,r},e_{p,t})=\hL.
\end{equation}
For $k\in\Z$ we set
$m:=\frac {k+1}{2},\ n:= \frac {k-1}{2}$ if $k$ is odd, and
$m:=\frac {k+2}{2},\ n:= \frac {k-2}{2}$ if $k$ is even.
In both cases we obtain
that $\g(e_{k,r},e_{p,t})=\hL$.
\begin{lemma}\label{L:vecrs}
For $r\ne t$ the value of the cocycle 
$\g(e_{k,r},e_{p,t})$ is given as a universal linear combination
of values of the cocycle at higher level.
\end{lemma}

Hence again only $r=s=t$ is of importance and we will drop the second
index. We obtain
\begin{equation}\label{E:motherv}
(m-n)\cdot \g(e_{n+m},e_p)+
(p-m)\cdot \g(e_{m+p},e_n)+
(n-p)\cdot \g(e_{n+p},e_m)=\hL.
\end{equation}
If we set $n=0$ and  use the antisymmetry we obtain
\begin{equation}
(m+p)\g(e_m,e_p)+(p-m)\g(e_{m+p},e_0)=\hL.
\end{equation}
Hence,
\begin{lemma}
If the level $l=m+p\ne 0$ then 
\begin{equation}
\g(e_m,e_p)=\frac {m-p}{m+p}\cdot\g(e_{m+p},e_0)+\hL.
\end{equation} 
\end{lemma}
It remains to deal with the level zero case.
Clearly $\g(e_0,e_0)=0$ due to antisymmetry.
Setting $p=-(n+1)$  and $m=1$  in \refE{motherv}
we get
\begin{equation}\label{E:vecrec}
(n-1)\g(e_{n+1},e_{-(n+1)})
=(n+2)\g(e_{n},e_{-n})
-(2n+1)\g(e_{1},e_{-1})+\hL.
\end{equation}
This recursion fixes the level zero starting from higher level and the
values of $\g(e_1,e_{-1})$ and  $\g(e_2,e_{-2})$.
\begin{proposition}\label{P:veczeron}
For a cocycle for the vector field algebra the level zero is given by the
data $\alpha_r$ and $\beta_r$ for $r=1,\ldots,K$ fixed by
\begin{equation}
\a_r=1/6\left(\g(e_{2,r},e_{-2,r})-2\g(e_{1,r},e_{-1,r})\right),\quad
\b_r=\g(e_{1,r},e_{-1,r})
\end{equation}
and higher level values via
\begin{equation}\label{E:veczeron}
\g(e_{k,r},e_{-k,s})=\left((k+1)k(k-1)\a_r+k\b_r\right)
\delta_r^s+\hL.
\end{equation}
Where $\hL$ denotes a universal linear combination of values of the cocycle
evaluated at  higher levels.

(b) A  cocycle which is bounded from above is uniquely given 
by the collection of values
\begin{equation}\label{E:valvec}
\g(e_{1,r},e_{-1,r}),\g(e_{2,r},e_{-2,r}), \ r=1,\ldots, K,
\quad
\g(e_{n,r},e_{0,r}),\ n\in\Z\setminus\{0\},\ r=1,\ldots, K.
\end{equation}
\end{proposition}
\begin{proof}
(a) For $r\ne s$  \refL{vecrs} gives the claim. For $r=s$ 
Equation \refE{vecrec} 
gives the recursive relation.
It remains to show the explicit formula. 
By antisymmetry it is enough to consider $k>0$ and there
it follows from induction
starting with $k=1$ and $k=2$.

Part (b) follows 
with the same arguments as in the proof of
\refP{zerom}.
\end{proof}
\begin{theorem}\label{T:vbound}
(a) Let $\g$ be a cocycle for the vector field algebra which is bounded from
above by $M$, then there exists a 
sum $W=\sum_{-\infty}^{n=M}\sum_r\beta_{n,r}\Omega^{n,r}$ of
quadratic differentials and a collection $\a_k\in\C$, $k=1,\ldots, K$ such that
$\g$ is the linear combination
\begin{equation}
\g(e,f)=\sum_{i=1}^K\a_i\g_i(e,f)+D_W(e,f)
\end{equation}
with 
\begin{equation}
D_W(e,f)=\sum_{-\infty}^{n=M}\sum_r\beta_{n,r}
\cint{C_i}\Omega^{n,r} [e,f]=\kndual {W}{[e,f]} .
\end{equation}
(b) The cohomology space of mixing cocycles bounded from above
is $K$-dimensional and generated by the classes
$[\g_r],\ r=1,\ldots,K$.
\end{theorem}
\begin{proof}
The proof of (a) is completely analogous the proof presented for the
mixing cocycles.
It allows to calculate $\beta_{n,r}$ recursively from above
to obtain any $\g(e_{n,r,}e_{0,r})$.
\newline
On level zero we calculate
\begin{equation}
\g(e_{1,r},e_{-1,r})=-2\beta_{0,r}+\hDT,
\quad 
\g(e_{2,r},e_{-2,r})=\frac 12\a_r -4\beta_{0,r}+\hDT.
\end{equation}
Hence,
\begin{equation}
\a_r=2\g(e_{2,r},e_{-2,r})-4\g(e_{1,r},e_{-1,r})+\hDT,\quad
\beta_{0,r}=-1/2\g(e_{1,r},e_{-1,r})+\hDT
\end{equation}
will realize the given values.
The argument for negative level is again the same as for the
mixing cocycle.
\newline
(b) follows again from (a) by \refP{indep}.
\end{proof}
%%%%%%%%%%%%%%%%%%%%%%%%%%%%%%%%%%%%%%%%%%%%%%%%%%%%%%
\begin{proof}[Proof of \refT{vector}.]
The proof has completely the same structure as the proof of the
function case and the mixing case respectively.
For the testing vector fields $f_n$,
and  $e_n$,  for infinitely many $n$ with $n\gg 0$ we require  the orders
$\ord_{P_k}(f_n)=n+1$,
$\ord_{P_1}(f_n)=-n-3g+1$,
$\ord_{P_k}(e_n)=-n+1$,
$\ord_{P_1}(e_n)=n-3g+1$
and the condition that they should be regular elsewhere and 
normalized at $P_k$.
Then 
\begin{equation}
\begin{gathered}
(\a_k-\a_1)\frac 1{24}(n+1)n(n-1)+A(e_n,f_n)=0,
\\
A(e_n,f_n):=D_V(e_n,f_n)-D_{V^*}(e_n,f_n).
\end{gathered}
\end{equation}
As above (see also the  appendix)
\begin{claim}\label{C:claimvec}
$A(e_n,f_n)=O(n)$.
\end{claim}
By letting $n$ go to $\infty$ we conclude that $\a_k=\a_1$ and 
obtain all the other results in the same way as for the mixing
cocycle.
In particular for $\alpha\ne 0$ we can suitable rescale $W$ and
incorporate it into the projective connection. This shows (b)
The behaviour of level zero \refE{veczero} follows from
\refP{veczeron} and the fact that for the separating cocycle
$\g_S$ we have $\g_S(e_{2,r},e_{-2,r})=1/2$ and  $\g_S(e_{1,r},e_{-1,r})=0$.
Part (d) follows from (a) with \refP{indep}.
\end{proof}

%%%%%%%%%%%%%%%%%%%%%%%%%%%%%%%%%%%%%%%%%%%%%%%%%%
\section{An application:
$\glib$ and wedge representations
for the  differential operator algebra}
\label{S:wedge}
%\input wedge.tex
%\section{Application:
%Local Cocycles for the  differential operator algebra
%and the 
%Wedge Representation}
%\label{S:wedge}
%\input wedge.tex
%%%%%%%%%%%%%%%%%%%%%%%%%%%%%

\subsection{The infinite matrix algebra $\glib$} $ $
\medskip

First let me recall the following facts about infinite-dimensional
matrix algebras (see \cite{KaRa} for details).
Let $\mathrm{mat}(\infty)$ be the vector space of (both-sided) infinite 
complex matrices. An element $A\in\mathrm{mat}(\infty)$
can be given as
\begin{equation}
A=(a_{ij})_{i,j\in\Z},\qquad a_{ij}\in\C.
\end{equation}
Consider the subspaces
\begin{equation}
\begin{aligned}
\gli&:=\{A=(a_{ij})\mid  \exists r=r(A)\in\N: a_{ij}=0, \text{ if } |i|,|j|>r\},
\\
\glib&:= \{A=(a_{ij})\mid \exists r=r(A)\in\N: a_{ij}=0, \text{ if } 
|i-j|>r\}.
\end{aligned}
\end{equation}
The matrices in $\gli$ have ``finite support'', the matrices in $\glib$
have  ``finitely many diagonals''.
The elementary matrices $E_{kl}$ are given as
$$
E_{kl}=(\de_i^k\de_j^l)_{i,j\in\Z}.
$$
For $\mu=(\ldots,\mu_{-1},\mu_{0},\mu_{1},\ldots)\in\C^\Z$ and $r\in\Z$
set
\begin{equation}
A_r(\mu):=\sum_{i\in\Z}\mu_iE_{i,i+r}
\end{equation}
to denote a  diagonal matrix where the diagonal is shifted
by $r$ positions to the right.
The elements $\{E_{kl}\}$ are a basis of $\gli$, the elements $\{A_r(\mu)\}$
are a generating set for $\glib$.
The subspaces $\gli$ and $\glib$  of $\mathrm{mat}(\infty)$  become associative
algebras with the usual matrix product.
To see that the multiplication is well-defined 
for  $\glib$ the fact that  every element has only finitely many diagonals
is of importance.
(Note that $\mathrm{mat}(\infty)$ itself is not an algebra.)
With the commutator they become infinite-dimensional Lie algebras.
In the terminology of Kac  and Raina the algebra $\glib$ is $\bar a_\infty$.

The Lie algebra $\glib$ admits a standard 2-cocycle,  \cite{Fucohom}.
For $A=(a_{ij})\in\glib$ set $\pi(A)=(\pi(A)_{ij})$ the matrix
defined by
\begin{equation}
\pi(A)_{ij}\ :=\quad
\begin{cases}
a_{ij},&i\ge 0,j\ge 0
\\
0,&\text{otherwise}.
\end{cases}
\end{equation}
The cocycle is  defined by
\begin{equation}\label{E:glcyc}
\a(A,B):=\tr\big(\pi([A,B])-[\pi(A),\pi(B)]\big).
\end{equation}
Note that the matrix expression under the trace has
finite support, hence the trace is well-defined.
Restricted to the subalgebra $\gli$ the cocycle vanishes.
The following proposition is well-known.
E. g. a proof can be found in \cite{Fucohom}.
\begin{proposition}
The bilinear form $\a$ defines a cocycle which is not
cohomologous to zero.
The continuous cohomology $\H^2_{cont}(\glib,\C)$   
is one-dimensional and generated by $\a$.
\end{proposition}
Let $\glih$ be the central extension defined via the cocycle class $[\a]$.
\begin{proposition}\label{P:acyc}
The cocycle $\a$ is a multiplicative cocycle, i.e.
\begin{equation}\label{E:acyc}
\alpha(A\cdot B,C)+
\alpha(B\cdot C,A)+
\alpha(C\cdot A,B)=0.
\end{equation}
\end{proposition}
\begin{proof}
Let us decompose the matrices $A$, $B$ and $C$ into the following
four boxes
\begin{equation}
X=\begin{pmatrix}
X_1&X_2
\\  X_3&X_4
\end{pmatrix},\qquad
\text{with}\quad
X_4=\pi(X).
\end{equation}
We will not distinguish between the boxes and the matrices 
in $\glib$ obtained by filling them up again to elements of $\glib$.
In particular $X=X_1+X_2+X_3+X_4$.
The matrices $X_2$ and $X_3$  have finite support.
A direct calculation shows
\begin{equation}
\a(X,Y)=\tr(X_3Y_2)-\tr(Y_3X_2).
\end{equation}
Hence
\begin{equation}\label{E:alaus}
\a(AB,C)=
\tr(A_3B_1C_2)+\tr(A_4B_3C_2)
-\tr(C_3A_1B_2)-\tr(C_3A_2B_4).
\end{equation}
Because all products have finite 
support all the traces make sense.
Permuting \refE{alaus} cyclically and adding the results  
gives \refE{acyc}.
\end{proof}
For the generators of $\glib$ we calculate
\begin{equation}\label{E:ars}
\a(A_r(\mu),A_{-s}(\mu'))=0,\quad r\ne s,
\end{equation}
and
\begin{equation}\label{E:ar}
\a(A_r(\mu),A_{-r}(\mu'))
\ =\quad
\begin{cases}
\qquad \quad 0,&r=0
\\
\quad \ \ \sum_{k=0}^{-1-r}\mu_k\mu'_{k+r},&r<0
\\
\quad-\sum_{k=0}^{r-1}\mu_{k-r}\mu'_k,&r>0.
\end{cases}
\end{equation}
For the basis elements of $\Fl$ we introduce a linear order
in a lexicographical way, i.e. $(n,r)>(m,s)$ if $n>m$ or ($n=m$ and 
$r>s$). Set $v_{Kn+r}:=f^{\l}_{n,r}$.
In this way we can assign by the $\A$, $\L$ or $\Do$-module structure
of $\Fl$ to every element of $A$ and $\L$ a infinite matrix in the usual
way if we use the basis elements $v_j$ together with its numbering.
The almost-grading of the module structure guarantees that the matrix
will be in $\glib$.
Denote the induced Lie homomorphism or the homomorphism 
of associative algebras by $\Phi_\l$.
By the almost-graded structure we can write
\begin{equation}\label{E:embedd}
\Phi_\l(A_{n,r})=\sum_{r=-K(n+L_1)}^{-Kn}A_r(\mu),
\qquad
\Phi_\l(e_{n,r})=\sum_{r=-K(n+L_2)}^{-Kn}A_r(\mu'),
\end{equation}
with  elements $\mu,\mu'\in\C^\Z$ given by the structure constants.
The numbers $L_1$ and $L_2$ are the upper bounds for the almost-graded
structure.
The cocycle $\alpha$ can be pulled back to $\A,\L$ and $\Do$
to obtain a cocycle $\gamma_\l$  by
\begin{equation}
\g_\l(e,g)=\Phi_\l^*(\a)(e,g)=\a(\Phi_\l(e),\Phi_\l(g)).
\end{equation}
\begin{proposition}
The cocycle $\g_\l$ obtained by pulling back $\a$ is a local cocycle which
is bounded from above by zero.
As cocycle of $\A$ it is multiplicative.
\end{proposition}
\begin{proof}
By \refE{embedd},\refE{ars} and \refE{ar} we see that it is indeed local and
bounded by zero from above.
\refP{acyc} shows that it is multiplicative.
\end{proof}
Let me remark that the multiplicativity follows also indirectly
because $\g_\l$ on $\A$ is obtained by restriction of a differential
cocycle and is local, see \refT{fml}.
But the property expressed in \refP{acyc}
is  also important  in more general situation, \cite{Shsc}.
\begin{theorem}\label{T:pullcyc}
The cocycle $\g_\l=\Phi^*_\l(\a)$ can be written as 
the following 
linear combination of the separating cocycles introduced above,
\begin{equation}\label{E:pullcyc}
\g_\l=\Phi^*_\l(\a)=
-\left(\g_S^{(f)}+\frac {1-2\l}{2}\g_{S,T_\l}^{(m)}+
2(6\l^2-6\l+1)\g_{S,R_\l}^{(v)}\right),
\end{equation}
with a suitable meromorphic affine connection $T_\l$ and a projective
connection $R_\l$ without poles outside of $A$ and at most
poles of order one at the points in $I$ for $T_\l$ and order two for
 $R_\l$.
\end{theorem}
\begin{proof}
The existence of such a linear combination 
with  possible coboundary terms follows from the
uniqueness results of  \refS{results}.
It remains to calculate the scalar factors.
But from the explicit expressions \refE{ar} of the cocycle $\alpha$
we calculate
immediately
\begin{gather*}
\g_\l(A_{1,r},A_{-1,r})=1,\\
\g_\l(e_{1,r},e_{-1,r})=-\l(\l-1),\quad
\g_\l(e_{2,r},e_{-2,r})=-(1-2\l)^2+2\l(2-2\l),\\
\g_\l(e_{1,r},A_{-1,r})=\l-1,\quad
\g_\l(e_{-1,r},A_{1,r})=\l.
\end{gather*}
The only structure constants necessary for the above calculations are
the values given in  \refP{boundary}.
We use \refE{funzero},\refE{veczero} and \refE{mixzero}
to calculate the factors in the combination.
All factors in front of the basic separating cocycles are
non-zero and the coboundary terms can be incorporated into
the connections $T_\l$ and $R_\l$.
Note that the overall minus sign could be removed by rescaling
the central element $t$ by $(-1)$.
\end{proof}
\begin{remark}
Recall that the three separating cocycles in \refE{pullcyc}
are linearly independent. Hence, the central extensions
$\widehat{\mathcal{D}}_\lambda^1$ of the differential operator
algebra  associated
to different weights $\l$ are not even after rescaling
the central element equivalent.
If we consider only the centrally extended $\A$ we see that the
same central extension $\widehat{\mathcal{A}}$ will do.
Clearly, the obtained central extensions 
$\Lh_\lambda$ to different $\l$ 
of $\L$ will be after rescaling of the central element
be equivalent.
But the explicit element in the class will depend on the 
weight $\l$ via the projective connection $R_\l$.
\end{remark}
\begin{remark}
In this article we considered only $\l\in\Z$.
But for $\l\in \frac 12\Z$ the formula \refE{pullcyc} will also be true
with the only exception of $\l=1/2$.
Here the mixing cocycle will vanish, but a boundary term 
$E_{V_{1/2}}$ will remain. Hence
\begin{equation}
\g_{1/2}=
-\g_S^{(f)}+\g_{S,R_{1/2}}^{(v)}+E_{V_{1/2}}.
\end{equation}
\end{remark}

Let me indicate the relevance of \refT{pullcyc}.
In quantum field theory one is usually searching for highest weight
representations of the symmetry algebra.  The modules $\Fl$ are clearly not 
of this type.
But there is procedure (which for the classical situation is
well-known) how to construct from $\Fl$ the space of semi-infinite
wedge forms $\Hl$ and to extend the action to it.
The naively extended action will not be well-defined.
It has to be regularized.
See \cite{SchlDiss} and \cite{Schlwed} for the details.
As in the classical case the regularization is done by
embedding the algebras via the action on an ordered basis
of $\Fl$ into $\glib$ and by using the standard regularization
procedure there. One obtains for $\glib$ only a projective
action which can be described as a linear representation of
the centrally extended algebra $\glih$ defined via the cocycle $\a$.
Pulling back the cocycle we obtain an action of $\Ah,\Lh$ and $\Dh$
on $\Hl$.
We are exactly in the situation discussed above.
\begin{theorem}
The space of semi-infinite wedge forms $\Hl$ carries a representation
of centrally extended algebras $\Ah$, $\Lh_\l$ 
$\widehat{\mathcal{D}}_\lambda^1$.
The defining cocycle $\g_\l$ for the central extension is given by 
\refE{pullcyc}.
The cocycles for  $\Ah$ and $\Lh_\l$ are  obtained by restricting
$\g$ to the subalgebras.
\end{theorem}
%%%%%%%%%%%%%%%%%%%%%%%%%%%%%%%%%%%%%%%%%%%%%%%%%%%
In \cite{SchlDiss,Schlwed}
the algebra of differential operators  $\Dal$ on $\Fl$ of 
arbitrary degree is introduced.
The action can be extended to $\Hl$ if we pass to the
central extension obtained from pulling back $\a$.
Hence we obtain a central extension $\Dalh$ and a cocycle for this 
algebra.
\begin{proposition}
The algebra 
$\Dal$ of meromorphic differential operators holomorphic outside $A$
admits a central extension $\Dalh$. The restriction of the 
defining cocycle to the subalgebra $\Do$ of differential
operators of degree less or equal one is given by 
\refE{pullcyc}.
The algebra 
 $\Dalh$ can be realized as operators on the space of 
semi-infinite wedge forms.
\end{proposition}
Further details will appear in \cite{Schlwed}.
In the classical situation this extension is the extension
given by the Radul cocycle \cite{Radcoc}.
Note that Radul gave it only for $\l=0$.
Again for the classical
situation and again only for $\l=0$, 
 Li \cite{Li} showed  that this is the only linear combination of
cocycles for $\Do$ which can be extended to  $\mathcal{D}_0$.

%%%%%%%%%%%%%%%%%%%%%%%%%%%%%
\section{An application:
Cocycles for the affine algebra}
\label{S:affine}
%\input affine.tex
%\section{Application:  %Cocycles for the affine algebra}
%\label{S:affine}
%\input affine.tex
%%%%%%%%%%%%%%%%%%%%%%%%%%%%%%%%%%%%%%%%%%%%%%%%%%

Let $\ga$ be a reductive finite-dimensional Lie algebra with a fixed
invariant, symmetric bilinear form $(.,.)$,
i.e. a form obeying $([x,y],z)=(x,[y,z])$.
Further down we will assume nondegeneracy.
For the semi-simple case the  Cartan-Killing form will do.
The {\it multi-point higher genus current algebra} (or  {\it multi-point
higher genus
loop  algebra}) is defined as
\begin{equation}
\gb:=\ga\otimes \A,\quad
\text{with Lie product}\quad
[x\otimes f,y\otimes g]:=[x,y]\otimes f\cdot g\ .
\end{equation}
We introduce a grading in  $\gb$ by defining
\begin{equation}
\deg(x\otimes A_{n,p}):=n.
\end{equation}
This makes  $\gb$ to an almost-graded Lie algebra.

Important classes of central extensions $\gh$ of $\gb$ are given by
\begin{equation}\label{E:eaff}
[\widehat{x\otimes f},\widehat{y\otimes g}]=
\widehat{[x,y]\otimes (f g)}+(x,y)\cdot\gamma_C (f,g)\cdot t,\qquad
[\,t,\gh]=0\ ,
\end{equation}
where 
$
\gamma_C(f,g)=\cint{C} fdg
$ 
is a geometric cocycle for the function algebra obtained by
integration along the cycle $C$.
As usual I 
set $\widehat{x\otimes f}:=(0,x\otimes f)$.
These algebras are called the higher genus (multi-point)
affine Lie algebras (or \KN\ algebras of affine type).

In the classical situation these are nothing else then the
usual  affine Lie algebras (i.e. the untwisted affine Kac-Moody algebras). 
For higher genus such algebras were introduced by Sheinman 
\cite{Shea,Shaff} for the two
point situation and by the author
for the multi-point situation \cite{Schlct, Schlhab}.
See also Bremner \cite{Bremce, Brem4a} for related work.
From the purely algebraic context, i.e $\A$ an arbitrary
commutative algebra without a grading they were studied
earlier by  Kassel \cite{Kaskd}, Kassel and Loday \cite{KasLod},
and others.
For the $C^\infty$-case see also Pressley and Segal \cite{PS}.

From \refT{function}  we immediately get
\begin{proposition}
Assume that $(.,.)$ is nondegenerate, then 
the cocycle $(x,y)\cdot \g_C(f,g)$ is local if and only if the 
integration cycle $C$ is a separating cycle $C_S$.
\end{proposition}
We might even  assume a more general situation:
\begin{proposition}
Let $\ga$  be a finite-dimensional Lie algebra which fulfills the 
condition $[\ga,\ga]\ne 0$.
Let $\g$ be a local cocycle for the current algebra $\gb$ of the
$\ga$.
Assume that there is a nondegenerate
invariant symmetric bilinear form $(.,.)$ on $\ga$ and
a bilinear form $\g^{(f)}$ on $\A$ such that $\g$ can be written as
\begin{equation}
\g(x\otimes f,y\otimes z)=(x,y)\cdot\g^{(f)}(f,g),
\end{equation} 
then $\g^{(f)}$ is a multiple of the separating cocycle
for the function algebra.
\end{proposition}
\begin{proof}
First,  $\g^{(f)}$ is obviously antisymmetric and hence a cocycle for
$\A$.
We calculate 
\begin{equation}
\g([x\otimes f,y\otimes g],z\otimes h)=\g([x,y]\otimes f\cdot g,z\otimes h)
=([x,y],z)\,\g^{(f)}(f\cdot g,h).
\end{equation}
For the  cocycle condition for the elements
$x\otimes f,y\otimes g$ and $z\otimes h$
we have to permute this cyclically and add the result up.  We obtain
(using the invariance of $(.,.)$)
\begin{equation}
([x,y],z)\left(\g^{(f)}(f\cdot g,h)+
\g^{(f)}(g\cdot h,f)+\g^{(f)}(h\cdot f,g\right)=0.
\end{equation}
By the condition $[\ga,\ga]\ne 0$ and by the nondegeneracy
of $(.,.)$ it follows that $\g^{(f)}$ is a multiplicative cocycle.
Applying \refT{function} yields the claim.
\end{proof}
\begin{theorem}\label{T:simple}
Let $\ga$ be a finite-dimensional simple Lie algebra 
with Cartan-Killing form $(.,.)$, then every
local cocycle for the current algebra $\gb=\ga\otimes \A$ is
cohomologous to a cocycle given by
\begin{equation}\label{E:simple}
\gamma(x\otimes f,y\otimes g)=a\cdot \frac{(x,y)}{2\pi\i}\int_{C_S}fdg,\quad
\text{with}\ a\in\C.
\end{equation}
In particular,
$\H_{loc}^2(\gb,\C)$ is one-dimensional and up to equivalence 
and rescaling there is
only one nontrivial local central extension $\gh$ of $\gb$.
\end{theorem}
\begin{proof}
Kassel \cite{Kaskd} proved that the algebra $\gb=\ga\otimes\A$ for any
commutative algebra $\A$ over $\C$ and any $\ga$ a simple 
Lie algebra admits
a universal central extension. It is given by
\begin{equation}
\gh^{\, univ}=(\Omega^1_{\A}/d\A)\oplus\gb,
\end{equation}
with Lie structure
\begin{equation}
[x\otimes f, y\otimes g]=[x,y]\otimes fg+(x,y)\overline{fdg},\qquad
[\Omega^1_{\A}/d\A,\gh^{\, univ}]=0.
\end{equation}
Here $\Omega^1_{\A}/d\A$ denotes the vector space of K\"ahler differentials
of the algebra $\A$.
The elements in $\Omega^1_{\A}$ can be given as $fdg$ with $f,g\in\A$, 
and $\overline{fdg}$ denotes its class modulo $d\A$.
This universal extension is not necessarily one-dimensional.
Let $\gh$ be any one-dimensional central extension of $\gb$. It will  be given
as a quotient of $\gh^{\,univ}$. Up to equivalence it can be given by
a Lie homomorphism $\Phi$
\begin{equation}
\begin{CD}
\gh^{\,univ}=\Omega^1_{\A}/d\A\oplus\gb
@>\ \Phi=(\varphi,id)\ >>\gh=\C\oplus\gb
\end{CD}
\end{equation}
with a linear form $\varphi$ on $\Omega^1_{\A}/d\A$.
The structure of $\gh$ is then equal to
\begin{equation}
[x\otimes f, y\otimes g]=[x,y]\otimes fg+(x,y)\varphi(\overline{fdg}) t,\qquad
[t,\gh]=0.
\end{equation}

In our situation $M\setminus A$ is an affine curve and 
$\Omega^1_{\A}/d\A$ is the first cohomology group 
of the complex of
meromorphic functions on $M$ which are holomorphic on $M\setminus A$
(similar arguments can be found in an article by Bremner \cite{Bremce}).
By Grothendieck's algebraic deRham theorem \cite[p.453]{GH}
the cohomology of the complex is isomorphic to 
 the singular cohomology of $M\setminus A$.
Hence such a linear form $\varphi$ can be given by choosing a linear
combination of cycle classes in $M\setminus A$ and integrating
the differential $\overline{fdg}$ over this combination.
By \refT{function} the locality implies that the combination
is a multiple of the separating cocycle.
\end{proof}
%%%%%%%%%%%%%%%%%%%%%%%%%%%%%%%%%%%
\begin{remark}
There is a warning in order.
The claim of the above theorem is not true for $\ga$ only reductive.
As a nontrivial example take $\ga=gl(n)$ and $\psi$ any antisymmetric
bilinear form on $\A$. Then
\begin{equation*}
\gamma(x\otimes f,y\otimes g)=\tr(x)\tr(y)\psi(f,g)
\end{equation*}
defines a cocycle. But $\psi$ can be chosen to be local without being
a geometric cocycle.
\end{remark}
Further details will appear in a forthcoming paper
\cite{SchlHg}.

%%%%%%%%%%%%%%%%%%%%%%%%%%%%%%%%
\appendix
%\newpage
\label{S:appendix}
\section{Asymptotic expansions}
%\input appendix.tex
%%% Appendix
%%%%%%%%%%%
In this appendix 
I show \refCl{claimmix}. The proofs of Claims \ref{C:claimvec},
\ref{C:indm}, and \ref{C:indv} are completely analogous.
Recall the claim: $A(e_n,g_n)=O(n)$.
First we deal with the genus $g=0$ case.
We might assume that $P_1$ corresponds to $z=0$,
$P_k$ to $z=1$, and $Q_1$ to $z=\infty$.

For \refCl{claimmix} we have
$
g_n=z^{-n}(z-1)^{n},\quad
e_n=z^{n}(z-1)^{-n+1}\fpz.
$
This implies: $e_ng_n'=z^{-1}\cdot n$.
Hence the claim.
(For \refCl{claimvec} 
the elements are given by
$
f_n=z^{-n+1}(z-1)^{n+1}\fpz,\quad
e_n=z^{n+1}(z-1)^{-n+1}\fpz,
$
and we  calculate  $[e_n,f_n]=2z(-n+nz)\fpz$.)

For genus $g\ge 1$
we restrict the presentation to  the case of generic positions of 
the points in $A$. For a nongeneric position there
might appear 
an additional factor.
It will not depend on $n$ hence it can be ignored in the analysis.
Also we might assume that $n\gg 0$ to avoid special
prescription necessary for small $n$ in the case of weight 0 and weight 1.

In \cite{Schleg} for certain elements of $\Fl$ explicit expressions
were given.
The formulas there are valid if the required orders at the points
sum up to $(2\lambda-1)g-2\l$.
This is exactly the case for the elements  considered here.

The  building blocks are (see \cite{Schleg})
\begin{enumerate}
\item
the prime form $E(P,Q)$, which is a multivalued 
holomorphic form on $M\times M$ of weight $-1/2$ 
in
each argument. It will vanish only  only along the diagonal;
the zero will be of first order,
\item
the $\sigma$-differential which is a multivalued holomorphic form 
of weight $g/2$ without zeros,
\item
the well-known $\vartheta$-function on the Jacobian of $M$,   
\item
the Jacobi map $J$, which embeds  $M$ into its Jacobian,
\item
the Riemann vector $\Delta\in\C^g$ (see \cite[I, p.149]{Mumtheta}.
\end{enumerate}

First we deal with the mixing cocycle situation.
We abbreviate
\begin{align*}
S(g_n,P)&:=
J(P)-(n+g)J(P_1)+nJ(P_k)+\Delta,\quad
\\
S(e_n,P)&:=J(P)+(n-3g)J(P_1)-(n-1)nJ(P_k)+J(Q_1)+3\Delta.
\end{align*}
For the case of a mixing cocycle the elements are  given as
(see \cite[ Equation (18)]{Schleg})
\begin{align*}
g_n:&=\beta_1^{-1}E(P,P_1)^{-n-g}E(P,P_k)^n\sigma(P)^{-1}\vartheta
(S(g_n,P)),
\\
&\text{with}\quad \beta_1:=E(P_k,P_1)^{-n-g}\sigma(P_k)^{-1}\vartheta
(S(g_n,P_k))\in\C,
\\
e_n&:=\beta_2^{-1}
E(P,P_1)^{n-3g}E(P,P_k)^{-n+1}E(P,Q_1)\sigma(P)^{-3}\vartheta
(S(e_n,P)),
\\
&\text{with}\quad \beta_2:=
E(P_k,P_1)^{n-3g}E(P_k,Q_1)\sigma(P_k)^{-3}\vartheta
(S(e_n,P_k))\in\C.
\end{align*}
By the genericity $\vartheta(S(e_n,P_k))$ and $\vartheta(S(g_n,P_k))$
will not be zero.

We calculate $e_ng_n'$ where we take the derivative with 
respect to the 
local variable at the point $P$.
We obtain 
\begin{multline*}
e_ng_n'=
\frac {E(P,P_1)^{-4g-1}}{E(P_k,P_1)^{-4g}}
\frac {\sigma(P)^{-5}}{\sigma(P_k)^{-4}}
\frac {E(P,Q_1)}{E(P_k,Q_1)}\times
\\
\times
\bigg((-n-g) E(P,P_k)\sigma(P)\frac {\vartheta(S(e_n,P))
\vartheta(S(g_n,P))}  {\vartheta(S(e_n,P_k))
\vartheta(S(g_n,P_k))} E'(P,P_1)
\\
+ (n-1) E(P,P_1)\sigma(P)\frac {\vartheta(S(e_n,P))
\vartheta(S(g_n,P))}  {\vartheta(S(e_n,P_k))
\vartheta(S(g_n,P_k))} E'(P,P_k)
\\
+E(P,P_1)\sigma'(P)\frac {\vartheta(S(e_n,P))
\vartheta(S(g_n,P))}  {\vartheta(S(e_n,P_k))
\vartheta(S(g_n,P_k))} E(P,P_k)+
\\
+E(P,P_1)\sigma(P)\frac {\vartheta(S(e_n,P))
\vartheta'(S(g_n,P)) J'(P)}  {\vartheta(S(e_n,P_k))
\vartheta(S(g_n,P_k))} E(P,P_k).
\bigg)
\end{multline*}
The $n$-dependence is only due to the obvious multiplicative 
factors of the first two terms and the quotients of the theta functions.
But for the latter quotients we obtain
\begin{align}\label{E:f1}
\left(\frac {\vartheta(S(e_n,P))
\vartheta(S(g_n,P))}  {\vartheta(S(e_n,P_k))
\vartheta(S(g_n,P_k))}\right)^{(k)}&=O(1),\quad k\in\N_0,
\\ \label{E:f2}
\left(\frac {\vartheta(S(e_n,P))
\vartheta'(S(g_n,P))J'(P)}  {\vartheta(S(e_n,P_k))
\vartheta(S(g_n,P_k))}\right)^{(k)}&=O(n),\quad  k\in\N_0.
\end{align}
To see this we choose a fundamental region for the theta function.
As long as the arguments are inside this region the expression will be
bounded. We only have to consider the automorphy factor appearing due
to the fact that we have to move back the point to the fundamental region.
Because $S(e_n,P)$ and $S(g_n,P)$ change exactly in opposite direction
these will cancel. Let $\Pi$  be  the period matrix 
of the curve. For the translation with  $\w=m_1+\Pi m_2, \ m_1,m_2 \in\Z^g$
we get  for \refE{f1}  
\begin{equation*}\label{E:f11}
\frac {\vartheta(S(e_n,P)-w)
\vartheta(S(g_n,P)+w)}  {\vartheta(S(e_n,P_k)-w))
\vartheta(S(g_n,P_k)+w)}
=
\frac {\vartheta(S(e_n,P))
\vartheta(S(g_n,P))}  {\vartheta(S(e_n,P_k))
\vartheta(S(g_n,P_k))},
\end{equation*}
and 
for \refE{f2} 
\begin{multline}\label{E:f22}
\frac {\vartheta(S(e_n,P)-w)
\vartheta'(S(g_n,P)+w)J'(P)}  {\vartheta(S(e_n,P_k)-w))
\vartheta(S(g_n,P_k)+w)}
=
\\
(-2\pi\mathrm{i}{}^t m_2J'(P))
\frac {\vartheta(S(e_n,P))
\vartheta(S(g_n,P))}  {\vartheta(S(e_n,P_k))
\vartheta(S(g_n,P_k))}
+
\frac {\vartheta(S(e_n,P))
\vartheta'(S(g_n,P))J'(P)}  {\vartheta(S(e_n,P_k))
\vartheta(S(g_n,P_k))}.
\end{multline}
In the last equation the term $m_2$ is responsible for the $O(n)$
behaviour.
This shows \refE{f1} and \refE{f2} for $k=0$.
Taking derivatives will not change the asymptotics because the derivatives
will always be symmetric with respect to $S(e_n,P)$ and  $S(g_n,P)$.
Hence these automorphy factors which could introduce higher powers of
$n$, i.e. the terms of the first kind of the r.h.s. of \refE{f22},
will cancel.
In the calculations 
 it is of importance that 
$S(e_n,P)-S(g_n,P)-S(e_n,P_k)+S(g_n,P_k)=0$.
Altogether this implies that also the derivatives are of the
required order.
This shows the claim.

%%%%%%%%%%%%%%%%%%%%%%%%%%%%%%%%%%%%%%%%%%%%%%%%%%
%
%   further sections:
%%%%%%%%%%%%%%%%%%%%%%%%%%%%%%%%%%%%%%%
%%%%%%%%%%  bibliography  %%%%%%%%%%
%\bibliographystyle{amsplain}
%\bibliography{/home/schlich/tex/bibtex/database}

%\end{thebibliography}
%\input kozykmain.bbl
%\providecommand{\bysame}{\leavevmode\hbox to3em{\hrulefill}\thinspace}

%%%%%%%%%%%%%%%%%%%%%%%%%%%%%%%%%%%%%%%%%%%%%%%%%%%%
\end{document}